\documentclass[number, 3p, times]{elsarticle}

\journal{Energy Economics}

\usepackage{amsmath}
\usepackage{amsthm}
\usepackage{amssymb}
\usepackage{mathtools}

\usepackage[table, dvipsnames]{xcolor}
\usepackage{booktabs}
\usepackage{chemformula}
\usepackage[ruled]{algorithm2e}
\usepackage{xurl}
\usepackage[hidelinks]{hyperref}
\usepackage{caption}
\usepackage{subcaption}
\usepackage[capitalise, noabbrev]{cleveref}

\crefformat{footnote}{#2\footnotemark[#1]#3}

\usepackage[UKenglish]{babel}
\usepackage{csquotes}
\usepackage[shortcuts]{extdash}
\usepackage[useregional]{datetime2}

\usepackage{graphicx}
\usepackage{graphbox}

\usepackage{parskip}
\usepackage{microtype}

\usepackage[section]{placeins}
\makeatletter
\AtBeginDocument{%
	\expandafter\renewcommand\expandafter\subsubsection\expandafter{%
		\expandafter\@fb@secFB\subsection
	}%
}
\makeatother

\DeclareMathOperator{\argmax}{argmax}
\DeclareMathOperator{\argmin}{argmin}

\newcommand{\Feps}{\mathcal{F}_{\varepsilon}}
\newcommand{\Fepsdash}{\mathcal{F}_{\varepsilon}'}
\newcommand{\Aeps}{\mathcal{A}_{\varepsilon}}
\newcommand{\Ninv}{N_{\text{inv}}}
\newcommand{\Nop}{N_{\text{op}}}
\newcommand{\copt}{c_{\text{opt}}}
\newcommand{\xrob}{x_{\text{rob}}}
\newcommand{\ych}{y_{\text{ch}}}
\newcommand{\allyears}{{\mathcal{Y}}}

\makeatletter
\def\ps@pprintTitle{%
    \let\@oddhead\@empty
    \let\@evenhead\@empty
    \def\@oddfoot{\footnotesize\itshape
        {Accepted version, submitted to Energy Economics. \copyright{} 2023 %
        \hspace{2ex} \href{https://creativecommons.org/licenses/by-nc-nd/4.0}{\includegraphics[align=c]{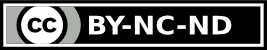}}}%
        \hfill\today}%
    \let\@evenfoot\@oddfoot
    }
\makeatother

\begin{document}

\begin{frontmatter}
	\title{Intersecting near-optimal spaces: European power systems with more resilience to weather variability}

	\author[1]{Aleksander Grochowicz\corref{cor1}\fnref{fn1}}
	\ead{aleksgro@math.uio.no}

	\author[2]{Koen van Greevenbroek\fnref{fn1}}

	\author[1]{Fred Espen Benth}
	\author[3]{Marianne Zeyringer}

	\cortext[cor1]{Corresponding author}
	\fntext[fn1]{Equal contribution}

	\affiliation[1]{%
		organization={Department of Mathematics, University of Oslo},
		addressline={P.O. Box 1053 Blindern},
		city={Oslo},
		postcode={0316},
		country={Norway}}
	\affiliation[2]{%
		organization={Department of Computer Science, UiT The Arctic University of Norway},
		addressline={Postboks 6050 Langnes},
		postcode={9037},
		city={Tromsø},
		country={Norway}}
	\affiliation[3]{%
		organization={Department of Technology Systems, University of Oslo},
		addressline={P.O. Box 70},
		postcode={2027},
		city={Kjeller},
		country={Norway}}

	\begin{abstract}
		\setlength\parindent{0pt}%
		\setlength\parskip{0.5\baselineskip}%
		\noindent We suggest a new methodology for designing robust energy systems.
		For this, we investigate so-called near-optimal solutions to energy system optimisation models; solutions whose objective values deviate only marginally from the optimum.
		Using a refined method for obtaining explicit geometric descriptions of these near-optimal feasible spaces, we find designs that are as robust as possible to perturbations.
		This contributes to the ongoing debate on how to define and work with robustness in energy systems modelling.

		We apply our methods in an investigation using multiple decades of weather data.
		For the first time, we run a capacity expansion model of the European power system (one node per country) with a three-hourly temporal resolution and 41 years of weather data.
		While an optimisation with 41 weather years is at the limits of computational feasibility, we use the near-optimal feasible spaces of single years to gain an understanding of the design space over the full time period.
		Specifically, we intersect all near-optimal feasible spaces for the individual years in order to get designs that are likely to be feasible over the entire time period.
		We find significant potential for investment flexibility, and verify the feasibility of these designs by simulating the resulting dispatch problem with four decades of weather data.
		They are characterised by a shift towards more onshore wind and solar power, while emitting more than 50\% less \ch{CO2} than a cost-optimal solution over that period.

		Our work builds on recent developments in the field, including techniques such as Modelling to Generate Alternatives (MGA) and Modelling All Alternatives (MAA), and provides new insights into the geometry of near-optimal feasible spaces and the importance of multi-decade weather variability for energy systems design.
		We also provide an effective way of working with a multi-decade time frame in a highly parallelised manner.
		Our implementation is open-sourced, adaptable and is based on PyPSA-Eur.
	\end{abstract}

	\begin{graphicalabstract}
		\includegraphics[scale=0.7]{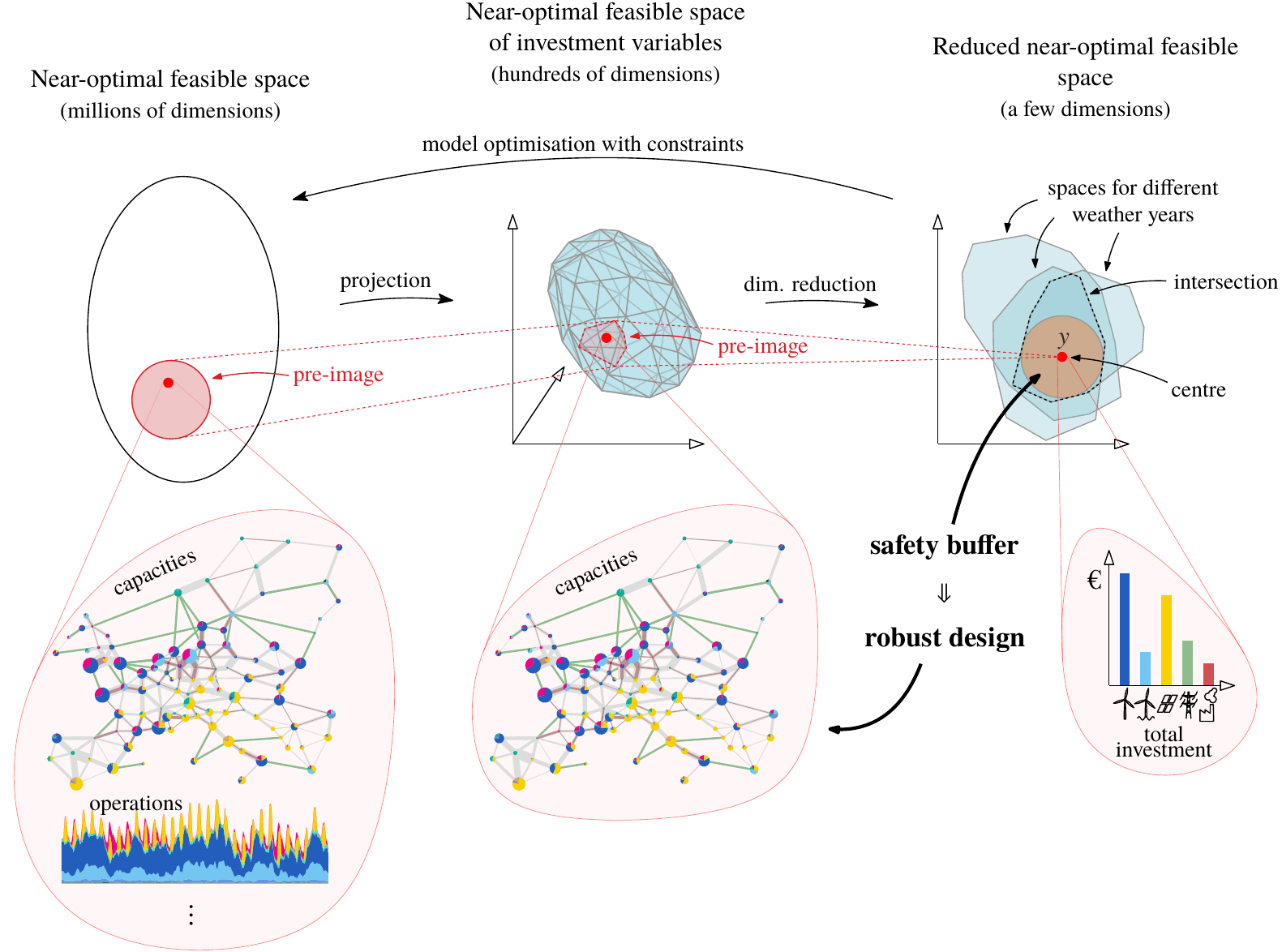}
	\end{graphicalabstract}
	
	\begin{highlights}
	\item Open energy system optimization model for Europe run with decades of weather data
	\item Different weather years impact optimal solutions, but near-optimal spaces intersect
	\item Geometry of near-optimal spaces allows for robustness against perturbations
    \item More investment lowers variable costs \& leads to more robustness and fewer emissions
    \item Approximating near-optimal spaces is parallelizable, making large problems tractable
	\end{highlights}
	
\end{frontmatter}

\section{Introduction}
\label{sec:introduction}

The climate crisis and tumbling prices for renewable technologies in the last decade are leading to an unprecedented shift to variable renewable and other low-carbon energy sources.
The pace and extent at which this transition is projected to take place often corresponds to a complete overhaul of currently existing energy systems within a few decades.
This necessitates a renewed understanding of the workings and planning of energy systems, taking into account future unknowns including weather, climate, costs, and politics.

In the domain of energy system modelling, grappling with these unknowns is an open research problem which has prompted many different approaches.
In particular, weather variability has been identified to have a large impact on model solutions: cost-optimal solutions are often fragile in the sense that modelling with one weather year can produce solutions that are infeasible for other weather years.
At the same time, recent work involving the relaxation of cost-optimality has revealed the opportunities and insights provided by the near-optimal feasible space of energy system optimisation models.
Still, existing methods only map out near-optimal solutions partially or heuristically, and a complete understanding of the theoretical and computational trade-offs has not yet been developed.

In this paper, we introduce a methodology to study robustness of energy systems against uncertainties in inputs and apply it to inter-annual weather variability.
To start with, we investigate methods for approximating projections of near-optimal feasible spaces of energy system models.
We then intersect these spaces for varying input data to produce solutions which are feasible in all scenarios.
Finally, we propose to consider points well in the \emph{interior} of this intersection as candidate robust solutions.
This geometric approach lends itself to building increased resilience against uncertainties and perturbation in a systematic way and generating greater flexibility for policymakers.

We are interested in the problem of long-term planning of the European power system with a large share of renewable production capacity in order to meet emission targets.
Capacity expansion models are often used for this purpose, in which capacities of generation technologies, storage, and transmission are optimised.
The model at the same time ensures feasibility by simultaneously solving the corresponding optimal dispatch problem over a certain time period, ``simulating'' the operations of the network.
Such models are vital exploratory tools used to shape high-level European energy policy.

In this study, we use the open-source, bottom-up energy system optimisation model (ESOM) called PyPSA-Eur \cite{PyPSAEur} for our implementation.
PyPSA-Eur consists of a model-building routine based on the PyPSA (Python for Power System Analysis) framework \cite{brown-horsch-ea-2018}, collecting and processing the required input data from various sources.
It assembles a faithful representation of the European high-voltage transmission grid and existing generation capacities, and uses Atlite \cite{hofmann-hampp-ea-2021} to compute capacity factor time series for renewable energy sources (PV, wind, hydro), based on historical ERA5 reanalysis weather data.
We use PyPSA-Eur with a partial greenfield approach under perfect foresight, including existing transmission (expandable), hydropower and nuclear capacities (both non-extendable), but optimising renewables (onshore \& offshore wind power, solar power), gas turbines and storage from zero.
For the purpose of demonstrating and validating our methodology, we use a spatial resolution of one node per country, and a 3-hourly temporal resolution (without time aggregation).
However, our techniques can be applied also when using a significantly higher resolution and including more technologies and energy sectors.

Usually, these types of energy system models are optimised using one or a few historical weather years \cite{ringkjob-haugan-ea-2018}, or different weather years are used for sensitivity analyses \cite{lombardi-pickering-ea-2020}.
The issue of weather year variability has been addressed in the literature \cite{pfenninger-staffell-2016,pfenninger-2017,collins-deane-ea-2018, staffell-pfenninger-2018,hilbers-brayshaw-ea-2020,craig-wohland-ea-2022,ruhnau-qvist-2022} and has been identified as an important factor for ESOM outcomes \cite{zeyringer-price-ea-2018}.
In particular, using only a single weather year as input data for ESOMs can produce design solutions which are over-fitted to that year, and are not feasible in general \cite{bloomfield-brayshaw-ea-2016,zeyringer-price-ea-2018}.
However, most previous studies running ESOMs with decades of weather data have either focused on a single country such as Germany \cite{ruhnau-qvist-2022}, the UK \cite{pfenninger-2017,zeyringer-price-ea-2018}, and the US \cite{dowling-rinaldi-ea-2020} (a single-node model), or used only a European dispatch model to solve for operations, not capacity expansion \cite{collins-deane-ea-2018}.
One study on the impact of different climate scenarios applied a TIMES model for the European power system with a decades-long modelling horizon~\cite{simoes-amorim-ea-2021} in which the usage of representative time slices limited the ability to model storage and capture medium- and long-term effects.
Pickering, Lombardi and Pfenninger have recently used the sector-coupled Euro-Calliope model in a study of the European energy system at a 96-node two-hourly resolution \cite{pickering-lombardi-ea-2022}; the main results were generated over a single weather year but validated using 8 additional weather years.

To the authors' knowledge, this is the first paper to run a spatially resolved ESOM for the European power sector with multiple decades of weather data, without aggregating to time slices.
For this analysis, we use 41 years of ERA5 reanalysis data \cite{era5-data} for the European continent (from 1980 until 2020 inclusive).
While we aim to find system designs which are feasible for all weather years under consideration, we base our methods on optimisations with single weather years in order to reduce the computational burden.
However, we are still able to optimise our model with all 41 weather years in order to validate our approach.

The robust solutions that we are interested in are near-optimal feasible solutions, meaning that their costs do not go beyond a previously defined threshold.
These solutions are ``close to cost-optimal'' and leave room for alternative objectives and desirable qualities; the additional costs we accept lie below the 9--23\% deviation from cost optimality (due to political, social, or technical reasons) that have been observed in recent years in the UK \cite{trutnevyte-2016}.
Following the works by Neumann and Brown \cite{neumann-brown-2021} and Pedersen et al.\ \cite{pedersen-victoria-ea-2021} we exploit the geometric shape and properties of the near-optimal space defined by the ESOM.
Instead of studying the full-dimensional near-optimal feasible space, we study a projection onto 5 relevant dimensions representing total investments in certain technologies.

By varying the weather years as inputs, we then construct one (reduced) near-optimal feasible space for each year.
When we intersect these near-optimal feasible spaces, we obtain a space of solutions in which each point represents a set of total investment decisions which are feasible for every year under consideration.

As the most robust candidate in the intersection, we choose the point which lays in the middle, being as far away from being infeasible as possible.
This means that changes in total investment decisions (up to a certain point) still leave us in the near-optimal feasible space for every weather year.
We then map the total investments back to a full system design, and verify its feasibility by simulating its operations over the entire time period.
Note that our form of ``robustness'' is a geometric concept (laying in the middle of a near-optimal feasible space) and is only loosely connected to robust optimisation.

Apart from contributing to the discussion on energy system robustness, our methods also have implications for ESOM parallelisation.
The difficulty in parallelising linear program (LP) solvers has been highlighted as the main barrier preventing ESOMs in taking advantage of increasing computational power~\cite{kotzur-nolting-ea-2021}.
While there are efforts to address this problem at the level of LP solvers~\cite{rehfeldt-hobbie-ea-2022}, we work at the level of model formulation.
Finding solutions which are feasible for many weather years by studying the intersection of their respective near-optimal spaces can be an alternative to solving ESOMs with many weather years outright, which is computationally prohibitive.
Thus, our methods constitute a way of heuristically replacing one large (difficult to parallelise) optimisation by many optimisations with single weather years.

Finally, we formalise and significantly deepen the understanding of the geometry and approximation of near-optimal feasible spaces of ESOMs.
Previously, Pedersen et al.\ proposed a methodology for approximating near-optimal feasible spaces of ESOMs \cite{pedersen-victoria-ea-2021}, and used the results to study the density of certain system design properties under projections of the near-optimal feasible space.
Furthermore, Lombardi et al.\ mapped out the utilisation of renewable capacity, transmission capacity, and storage capacity of chosen near-optimal solutions, depending on different uncertainties, indicating overlaps between these \cite{lombardi-pickering-ea-2020}.
We detail the dimension reductions involved in working with near-optimal feasible spaces, and how to map back and forth between the different stages.
We then propose several variations on a general algorithm for approximating reduced near-optimal feasible spaces, and analyse their convergence characteristics.
The application of geometric descriptions of near-optimal spaces to studying different weather years is also novel.

In summary, our paper contributes to the literature on ESOMs in several ways.
We formalise a general framework for working with and intersecting near-optimal feasible spaces, which allows us to study uncertainties of different kinds.
We apply this framework to a first-of-its-kind study of robustness of highly renewable scenarios for the European power system to decades of weather data.
Beyond robustness, the methods also contribute to parallelisation of energy system optimisation models.

In \cref{sec:methods-formal-defs} we formalise the methodology and introduce the necessary steps to define ``robust'' energy system designs.
Afterwards, in \cref{sec:implementation}, we describe the adapted PyPSA-Eur model we use and our modelling set-up to obtain power systems resilient to 41 years of weather data.
In \cref{sec:results} we present our main findings on using intersections of near-optimal feasible spaces, features of robust solutions, and performance.
We discuss ramifications of our approach in \cref{sec:discussion}, before we conclude with \cref{sec:conclusion}.

\section{Methodology and formal definitions}
\label{sec:methods-formal-defs}
We introduce the methodology used in this paper and describe how we apply this to obtain energy system designs robust to weather variability:
in \cref{sec:near-optimality} we revisit the concept of near-optimality.
In \cref{sec:dim-reduction} we discuss how dimension reduction is necessary to describe the near-optimal feasible space in a computationally tractable manner, and we elaborate on how we approximate near-optimal feasible spaces.
Afterwards we introduce robustness as the geometric property of lying in the intersection of different near-optimal feasible spaces (\cref{sec:robustness}).
In \cref{sec:chebyshev-centre} we then justify our choice of the Chebyshev centre as our robust solution of choice, as its location implies maximal stability to perturbations.
Finally, in \cref{sec:going-back} we suggest different allocations that can translate a point in the low-dimensional intersection of the near-optimal feasible spaces to a spatially resolved (full-dimensional) energy system design.

\subsection{Near-optimality}
\label{sec:near-optimality}

Let a capacity expansion problem be given as the linear program
\begin{equation}
	\label{eq:lin-prog-def}
	\min c \cdot x \quad\text{such that}\quad Ax \leq b.
\end{equation}
Here, $x \in \mathbb{R}^N$ is a vector of decision variables, $c \in \mathbb{R}^N$ the coefficients of the objective function, and $A \in \mathbb{R}^M \times \mathbb{R}^N$ and $b \in \mathbb{R}^M$ give the set of $M$ linear constraints.
Let $\mathcal{F}$ be the feasible space of the linear program \cref{eq:lin-prog-def}, defined as
\begin{equation}
	\mathcal{F} := \{ x \mid Ax \leq b\}.
\end{equation}
Letting $x^* \in \mathcal{F}$ be an optimal solution with objective value $c \cdot x^* = \copt \in \mathbb{R}$, and $\varepsilon > 0$ a chosen slack level, we define the \textit{$\varepsilon$-near-optimal feasible space} as
\begin{equation}
	\label{eq:slack-condition}
	\Feps := \{ x \in \mathcal{F} \mid c \cdot x \leq (1 + \varepsilon) \cdot \copt\}.
\end{equation}
When $\varepsilon$ is clear from the context, we simply refer to $\Feps$ as the \emph{near-optimal} space.
For general linear programs, $\Feps$ is a convex polyhedron, and when $x$ is bounded (as is the case for energy system models), $\Feps$ is a convex polytope.
To work with $\Feps$ geometrically, we can solve the optimisation problem $\min d \cdot x \text{ s.t. } Ax \leq b \text{ and } c \cdot x \leq (1 + \varepsilon) \cdot c_{\text{opt}}$ for some objective $d$ in order to find a vertex or boundary point of $\Feps$.
In the context of ESOMs, this amounts to solving the energy system model once with an alternative objective function.

The definition of a near-optimal space is not new in the context of ESOMs, and previous work has explored the near-optimal space either through uniform sampling as \cite{pedersen-victoria-ea-2021}, maximally different solutions as \cite{decarolis-2011,price-keppo-2017}, or extreme points of the space \cite{neumann-brown-2021}.

\subsection{Dimension reduction}
\label{sec:dim-reduction}

\begin{figure}
	\centering
	\includegraphics{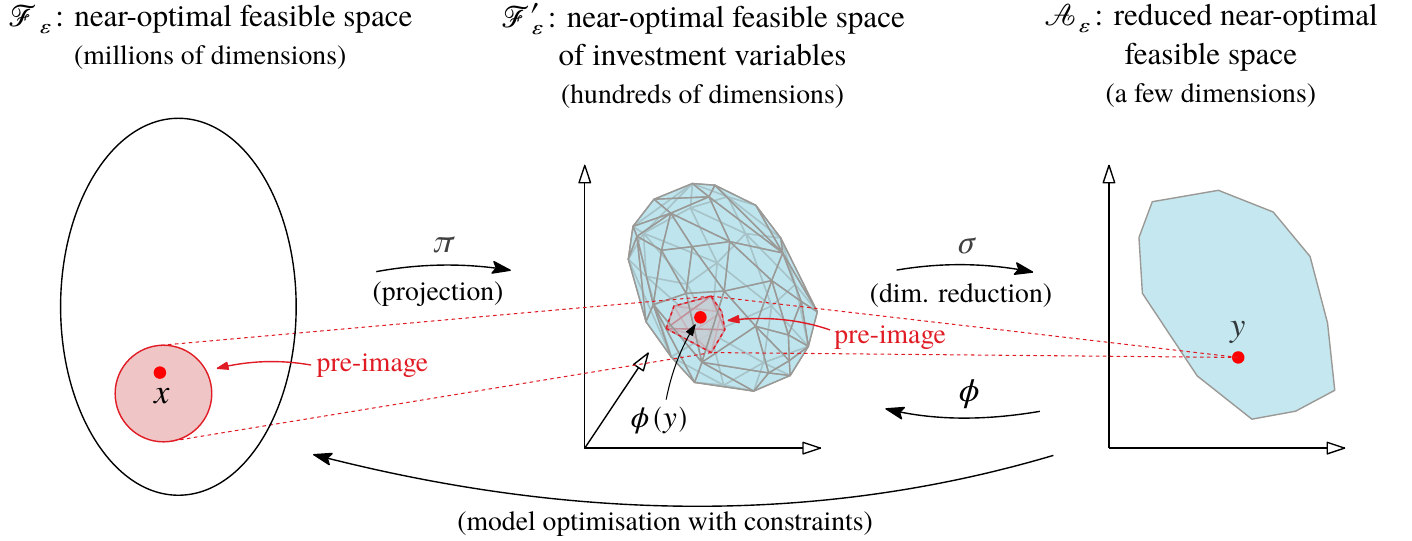}
	\caption{Illustration of how the spaces $\Feps$, $\Fepsdash$ and $\Aeps$ are connected.
		The map $\phi$ can be seen as a composition of a model optimisation ``$\min c \cdot x \text{ s.t. } Ax \leq b \text{ and } \sigma \circ \pi(x) = y$'' composed with the projection $\pi$. It gives an explicit model design based on chosen coordinates $y$ in the reduced near-optimal feasible space $\Aeps$.
	}
	\label{fig:dim-reductions}
\end{figure}

The near-optimal space $\Feps$ is high-dimensional and complex; large-scale ESOMs typically involve millions of variables (dimensions) and constraints (hyperplanes defining the polytope).
In this section we reduce to a much lower-dimensional space in two steps; see \cref{fig:dim-reductions} for an overview of the maps and spaces involved.

In the case of capacity expansion models, we are most interested in \emph{investment decision variables} of the linear program, as opposed to all other (operational) decision variables.
Specifically, let $x = (x^I,x^O)^{T} \in \mathbb{R}^N$ be split into investment decision variables $x^I \in \mathbb{R}^{\Ninv}$ and operational decision variables $x^O \in \mathbb{R}^{\Nop}$, where $N = \Ninv + \Nop$.
Then we define the projection map $\pi \colon \mathbb{R}^N \to \mathbb{R}^{\Ninv}$ by simply forgetting about the operational decision variables.
The image
\begin{equation}
	\label{eq:fepsdash-def}
	\Fepsdash = \pi(\Feps) = \{x^I \in \mathbb{R}^{\Ninv} \mid x = (x^I,x^O)^{T} \in \Feps\}
\end{equation}
of $\Feps$ under $\pi$ is the $\Ninv$-dimensional $\varepsilon$-near-optimal feasible space of investment variables.

The convex polytope $\Fepsdash$ consists of all points $x^I$ such that an energy system with capacity investments given by $x^I$ is feasible and whose total system cost (including operations) is at most $(1 + \varepsilon) \cdot \copt$.
In short, it is the space of all near-optimal feasible investment decisions.
This makes an explicit description of $\Fepsdash$ interesting for decision-makers in order to explore different kinds of near-optimal investments.

However, in a model with a high spatial resolution, the number of investment decision variables $\Ninv$ is typically still in the hundreds or more (with multiple investment decisions at each node, and transmission expansion).
This makes the polytope $\Fepsdash \subseteq \mathbb{R}^{\Ninv}$ difficult to work with, visually and mathematically.
Specifically, in order to work with $\Fepsdash$ we would want to find a set of points $P$ such that $\Fepsdash$ is the convex hull of $P$.
However, the number of vertices of an $\Ninv$-dimensional polytope defined by $M$ hyperplanes is in $O(M^{\lfloor \Ninv / 2 \rfloor})$ --- see~\cite{toth-goodman-ea-2017}, Chapter 26.
This puts a precise description of $\Fepsdash$ in terms of vertices out of reach.

One solution is to map down to a much lower-dimensional space where we group and aggregate investment decision variables.
Let $x^I = (x_1, \dots, x_{\Ninv})$ be the individual investment decision variables.
Let $T_1, \dots, T_k $ be a collection of sets of indices with $T_i \subseteq \{1, \dots, \Ninv\}$ and $T_i \cap T_j = \emptyset$.
For each index $j$ in one of these sets, we also choose a coefficient / weight $c_j$.
Then we define a linear map $\sigma \colon \Fepsdash \to \mathbb{R}^k$ as:
\begin{equation}
	\label{eq:aggregation}
	\sigma(x) = \textstyle \left(\sum_{j \in T_i} c_j x_j \right)_{i=1}^k.
\end{equation}

In our case, we take each $T_i$ to be the set of indices identifying decision variables that belong to a specific technology.
We weight these decision variables $x_j$ (for $j \in T_i$) by their respective capital costs $c_j$.
Specifically, throughout this paper we consider $k=5$, with $T_1, \dots, T_5$ corresponding to transmission expansion, PV expansion, onshore wind expansion, offshore wind expansion and gas turbine expansion respectively.
In effect, $\sigma$ maps a vector $x^I$ of investment decisions to a summary of selected total investment costs.

Let
\begin{equation}
	\label{eq:aeps-def}
	\Aeps = \sigma(\Fepsdash) = \{\sigma(x^I) \mid x^I \in \Fepsdash\} \subseteq \mathbb{R}^k
\end{equation}
be the image of $\Fepsdash$ under $\sigma$.
Then $\Aeps$ is a $k$-dimensional convex polytope (since convex polytopes are preserved by linear maps).
Note that in our specific choice of $T_1, \dots, T_5$ we have not included \emph{all} investment decision variables in $\sigma$, only those we deemed most important for the particular model instances we work with.
Of course, different dimension reductions can be achieved by other choices of aggregation (groups of indices $T_i$ and coefficients $c_{j}$).
While we have taken the coefficients $c_{j}$ to be capital costs (making investment in different technologies easier to compare), the coefficients could, for example, also be set to $1$ in order to consider only capacities.

The utility of the reduced near-optimal space $\Aeps$ is as a proxy for system feasibility.
If we can describe $\Aeps$ well, we can quickly assess whether any given set of total investments $y \in \Aeps$ can result in a feasible system design.
However, by aggregating investment decision variables, we lose information on the specific feasible system designs $\sigma^{-1}(y) \subseteq \Fepsdash$ realising the total investments $y$.
The trade-off is that the fewer dimensions $k$ we aggregate to, the easier $\Aeps$ is to work with, but the less information it gives us.
Each point $y \in \Aeps$ can have a large preimage under $\sigma$ and $\pi$, meaning there may be many near-optimal feasible solutions $x \in \Feps$ with the given total investments $y$ (see \cref{fig:dim-reductions}).
As we discuss in more detail in \cref{sec:going-back}, to find a specific solution $x \in (\sigma \circ \pi)^{-1}(y) \subseteq \Feps$, we must solve a version of the original model \cref{eq:lin-prog-def} in $\mathbb{R}^N$.

Similarly to \cite{pedersen-victoria-ea-2021}, we attempt to describe the low-dimensional space $\Aeps$ explicitly.
Again, in order to find an explicit description of the polytope $\Aeps$, we would like to find (a large subset of) its vertices.
We can do this by optimising over $\Aeps$ with different objective functions or \emph{directions} in $\mathbb{R}^k$.
For each direction $d \in \mathbb{R}^k$, the solution to the linear program
\begin{equation}
	\label{eq:new-objective}
	\max d \cdot y \quad\text{such that}\quad y = \sigma(\pi(x)) \text{ and } Ax \leq b \text{ and } x \in \Feps
\end{equation}
is an extreme point of $\Aeps$ in the direction $d$.
The above linear program can be solved by solving the original problem in \cref{eq:lin-prog-def} with the new objective function $- d \cdot \sigma(\pi(x))$ (which is linear), and mapping the solution to $\Aeps$ by $\sigma \circ \pi$.
In effect, for each extreme point of $\Aeps$ that we want to find, we need to solve the original capacity expansion problem once with an adapted objective function.

While we cannot expect to find all vertices of $\Aeps$ (quadratic in the number of model constraints when $k=4,5$ \cite{toth-goodman-ea-2017}), we want to find a set of extreme points $P$ such that their convex hull approximates $\Aeps$ well.
Given a ``budget'' of $n$ optimisations (and hence $n$ extreme points), the natural question is: how do we choose the directions $d_1, \dots, d_n$ to optimise in, in order to get a set of points $P$ whose convex hull approximates $\Aeps$ the best possible?

\begin{figure}
	\centering
	\includegraphics{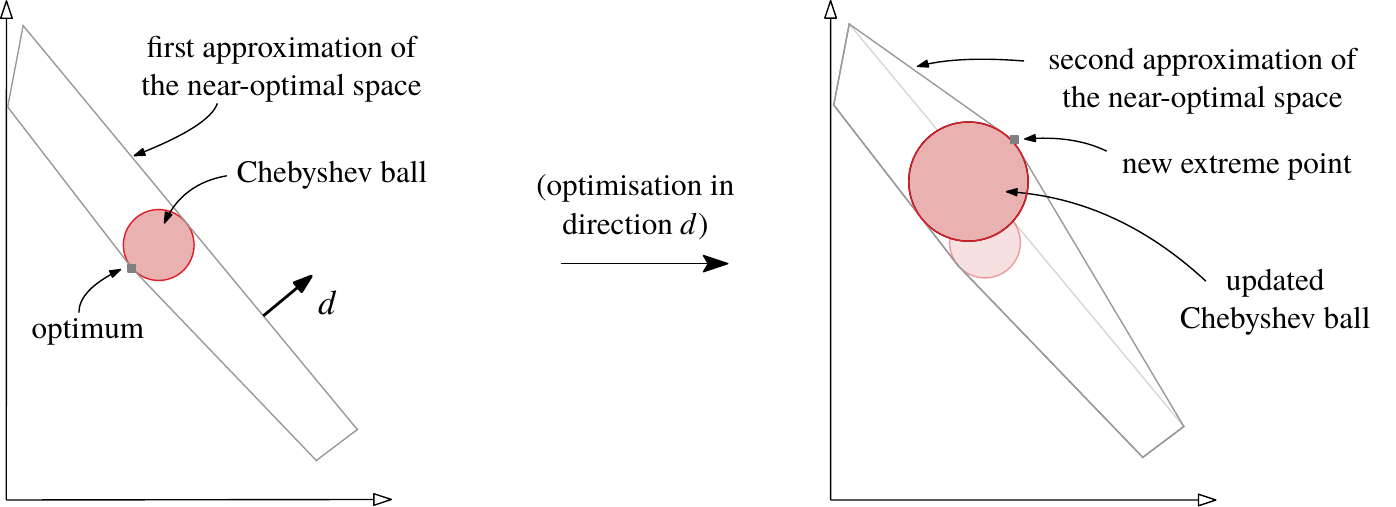}
	\caption{Illustration of a single step of \cref{alg:near-opt-approx}.
		Additionally, the Chebyshev ball (see \cref{sec:chebyshev-centre}) is shown at each stage.}
	\label{fig:directions-procedure}
\end{figure}

We use an iterative approach to approximate $\Aeps$ while filtering on already used (or similar) directions.
At each step we optimise in a different direction: depending on which property of the near-optimal feasible space is of interest, there can be many different ways to choose these directions.
In our case we are interested in the largest ball within the near-optimal feasible space, the Chebyshev ball (see \cref{sec:chebyshev-centre}).
Thus, we optimise at each iteration in the normal direction to a facet tangential to the Chebyshev ball of the current approximated polytope.
If these have been exhausted, we choose the normal direction to the largest facet by volume.
For other approaches to choosing directions, caveats and performance comparisons, see \cref{sec:direction-generation}.
An illustration of one step in this process is shown in \cref{fig:directions-procedure}.
A simplified version of the algorithm (based only on normals to large facets) is given in pseudo-code in \cref{alg:near-opt-approx}.

Once in the low-dimensional space $\mathbb{R}^k$ (with $k=5$ in our case), the complexity of the geometric objects of interest is low enough that we can do efficient exact computations on them.
In particular, recall that a polytope with $n$ vertices can only have $O(n^2)$ facets for $k=4, 5$, so computing the convex hull of the vertices, its volume, etc. has a time complexity of $O(n^2)$.
We use the qhull software \cite{barber-dobkin-ea-1996} for computational geometry related to polytopes.
Since the number of vertices $n$ which we can compute (by solving \cref{eq:new-objective} for each vertex) is limited, the complexity of $O(n^2)$ is acceptable for our purposes.
In particular, note that for our application, the complexity of \cref{alg:near-opt-approx} is dominated by the model optimisations, not the computational geometry.

\begin{algorithm}
	$P := \emptyset$\;
	\For{$d \in \{e_1, -e_1, e_2, -e_2, \dots, e_k, -e_k\}$}{
		Let $y$ be extreme point on $\Aeps$ in direction $d$ (optimisation)\;
		Add $y$ to $P$\;
	}
	\For{$i \in \{1, \dots, n\}$}{
		Let $H$ be the convex hull of $P$\;
		Let $F_1, F_2, \dots$ be the facets of $H$, sorted by decreasing volume\;
		Let $d_i$ be the normal of facet $F_i$\;
		Let $d$ be the first of $d_1, d_2, \dots$ which is not within a small angle $\theta$ of any previously used direction\;
		Let $y$ be extreme point on $\Aeps$ in direction $d$ (optimisation)\;
		Add $y$ to $P$\;
	}
	Return convex hull of $P$\;
	\caption{Outline of algorithm for approximating $\Aeps$}
	\label{alg:near-opt-approx}
\end{algorithm}

\subsection{Intersections and robust solutions}
\label{sec:robustness}

One of the new ideas we propose is to investigate the intersections of the near-optimal spaces of related capacity expansion problems, or different instances of the same abstract model.
Of course, if $\mathcal{F}^{(a)}, \mathcal{F}^{(b)} \subseteq \mathbb{R}^N$ are the feasible spaces of two linear programs $A$ and $B$, then $\mathcal{F}^{(a)} \cap \mathcal{F}^{(b)}$ is simply the space of all solutions $x$ which are feasible for both problems.
More interestingly for capacity expansion problems, consider ${\Fepsdash}^{(a)} \cap {\Fepsdash}^{(b)}$: the space of all investment allocations which are both feasible and near-optimal for both $A$ and $B$.

In our case, we consider the near-optimal spaces for optimisation problems defined with different weather years.
Specifically, we use 41 years of reanalysis weather data (1980-2020) in order to compute capacity factors and load time series as input for our model (\cref{sec:data} and \cref{sec:data-appendix}).
This gives 41 model instances, each defined with the input data from a different weather year.
For brevity, let $\allyears := \{1980, \dots, 2020\}$ denote the set of weather years.
Then for $i \in \allyears$ write
\begin{equation}
	\label{eq:lin-prog-def-years}
	\min c^{(i)} \cdot x^{(i)} \quad\text{such that}\quad A^{(i)} x^{(i)} \leq b^{(i)}
\end{equation}
for the LP in \cref{eq:lin-prog-def} defined with weather year $i$.
Let $\mathcal{F}^{(i)}$ be the feasible space of the above LP.
From these feasible spaces, we want to recover investment allocations which are feasible for each of the weather years.

Since total system costs vary considerably between cost optimisations with different weather years, we define a uniform system cost bound across all weather years.
Specifically, letting $c_{\text{opt}}^{(i)}$ be the optimal objective value (total system cost) for weather year $i$, we take
\begin{equation}
	\label{eq:uniform-slack}
	c^* = \max_{i \in \allyears} c_{\text{opt}}^{(i)}
\end{equation}
to be the highest optimal system cost across all weather years under investigation.
In defining near-optimal spaces with the different weather years, we then set the slack relative to $c^*$ instead of relative to each $c_{\text{opt}}^{(i)}$ for each weather year individually:
\begin{equation}
	{\Feps}^{(i)} = \{x \in \mathcal{F}^{(i)} \mid c \cdot x \leq (1 + \varepsilon) \cdot c^*\}.
\end{equation}
We define $\Fepsdash^{(i)}$ and $\Aeps^{(i)}$ similarly to $\Fepsdash$ and $\Aeps$ (\cref{eq:fepsdash-def,eq:aeps-def}).

Now, the intersection $\bigcap_{i  \in \allyears} \Fepsdash^{(i)}$ is the space of investment decisions which are feasible for all weather years under consideration and near-optimal (relative to the most expensive year)\footnote{In fact, this resembles the near-optimal feasible space of a robust optimisation program defined over the weather years $\allyears$. Strictly speaking, however, the different operational variables for different weather years make this a loose generalisation of classical robust optimisation.}.
We are interested in the intersection of the reduced near-optimal space $\Aeps^{(i)}$.
Repurposing our notation slightly, we write $\Aeps := \bigcap_i \Aeps^{(i)}$.
Note that $\Aeps = \bigcap_i \sigma(\Fepsdash^{(i)}) = \sigma\left(\bigcap_i \Fepsdash^{(i)}\right)$.
However, while we cannot easily find explicit descriptions of the spaces $\Fepsdash^{(i)}$, we \emph{can} approximate each (low-dimensional) $\Aeps^{(i)}$ as explained above.
This, in turn, enables us to find an explicit (approximate) description of $\Aeps = \bigcap_i \Aeps^{(i)}$.

We call a set of total investments $y \in \mathbb{R}^k$ \emph{robust} when $y \in \Aeps$.
This means that for each weather year $i$, there exists some near-optimal feasible model solution $x^{(i)} \in \Feps^{(i)}$ (including investment and operation decisions) such that $\sigma \circ \pi(x^{(i)}) = y$.
Note that the same total investment for each technology, $y \in \Aeps$, may be spread differently onto the different nodes of the model for each $x^{(i)}, i \in \allyears$.
Formally speaking, it is plausible that $\bigcap_i \Fepsdash^{(i)} = \emptyset$ even if we find some robust point $y \in \Aeps$.
In \cref{sec:going-back}, however, we propose different methods for finding robust allocations $x^{I} \in \bigcap_i \Fepsdash^{(i)}$ such that $\sigma(x^I) = y$ if they exist, and in \cref{sec:res-weather-years-feasibility} we show that this works well in practice.

\subsection{Chebyshev centre}
\label{sec:chebyshev-centre}
\begin{figure}
	\centering
	\begin{subfigure}[t]{0.48\textwidth}
		\centering
		\includegraphics{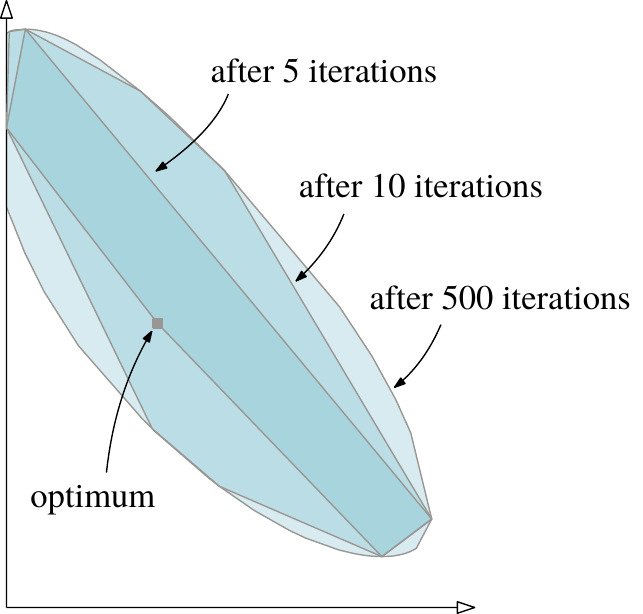}
		\caption{The approximation of a reduced near-optimal feasible space after 5, 10, and 500 iterations (with decreasing opaqueness).}%
		\label{fig:intersection-approx}
	\end{subfigure}
	\hfill
	\begin{subfigure}[t]{0.48\textwidth}
		\centering
		\includegraphics{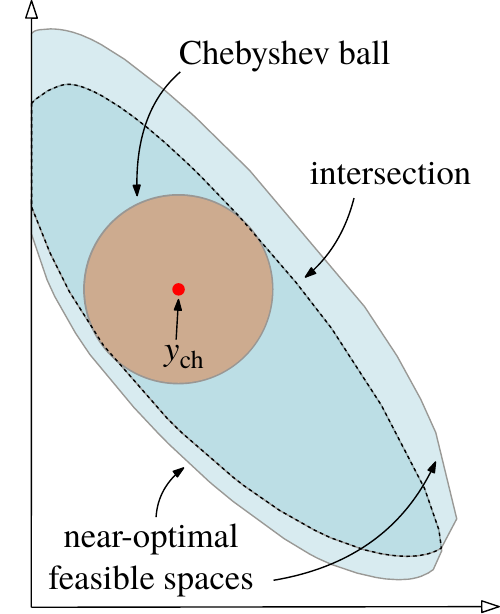}
		\caption{The Chebyshev centre $\ych$ of the intersection of two reduced near-optimal feasible spaces. The radius of the Chebyshev ball around the centre can give us an idea of how much flexibility we might expect around this point.}%
		\label{fig:intersection-chebyshev}
	\end{subfigure}
	\caption{Near-optimal space approximation, intersection and Chebyshev centre.}%
	\label{fig:approx-intersect-centre}
\end{figure}

Among the total investment decisions in the intersection $\Aeps = \bigcap_i \Aeps^{(i)}$, we want to find choices that are not only feasible for all years considered.
We want to find the most resilient choice among all the alternatives.
We therefore select the point $\ych \in \Aeps$ maximally removed in all directions from the boundary of $\Aeps$, meaning $\ych$ is as far away from being infeasible as possible.
This is realised if we pick $\ych$ to be the Chebyshev centre (see~\cite{boyd-vandenberghe-2004}, Section 8.5.1), i.e.
\begin{equation}
	\ych = \argmax_{y \in \Aeps} (r) \text{ s.t. } B_r(y) \subseteq \Aeps,
\end{equation}
where $B_r(y)$ is the ball of radius $r$ around $y$.
\cref{fig:directions-procedure,fig:intersection-chebyshev} show examples of Chebyshev balls.
The point $\ych$ can be found efficiently using a linear program.
Specifically, let $a_j$ be the normal vectors of the hyperplanes supporting $\Aeps$ and $b_j$ the associated offsets, so that each $y \in \Aeps$ satisfies $a_j \cdot y \leq b_j$ for all $j$.
Then $\ych$ is given by
\begin{equation}
	\label{eq:chebyshev-lp}
	\max r \quad\text{such that}\quad a_j \cdot y + r \| a_j \| \leq b_j \quad \forall j \;\text{ and }\; r \geq 0.
\end{equation}

\subsection{Disaggregating robust solutions}
\label{sec:going-back}

The previous steps leading to the Chebyshev centre, one chosen robust point, have all been performed in the reduced $\varepsilon$-near-optimal feasible spaces and their intersection, $\Aeps$.
The designs that are of ultimate interest to us, however, are elements in $\Fepsdash$, including \emph{all} investment decision variables.
Thus, we would like to define a function $\phi: \Aeps \to \Fepsdash$ mapping robust total investments to complete system designs realising those total investments.
Being more precise, we need to specify over which weather years the space we map back to is defined --- we want a map
\begin{equation}
	\label{eq:going-back}
	\phi: \Aeps \to \Fepsdash^\allyears,
\end{equation}
recalling the notation $\allyears = \{1980, \dots, 2020 \}$.
From this map, we want to obtain a robust near-optimal energy system design $\xrob = \phi(\ych)$ at the Chebyshev centre of $\Aeps$.

We can write our ESOM as in \cref{eq:lin-prog-def-years} but instead define it with all 41 weather years as $\min c^\allyears \cdot x^\allyears \text{ such that } A^\allyears x^\allyears \leq b^\allyears$; its $\varepsilon$-near-optimal feasible space of investment variables is $\Fepsdash^\allyears$.
Then we define the map $\phi$ (which we call $\phi^\allyears$ for clarity) by adding constraints to the above LP to ensure that the solution has total investments given by $y \in \Aeps$:
\begin{equation}
	\label{eq:exact-map}
	\phi^\allyears(y) = \pi(\argmin c^\allyears \cdot x^\allyears) \;\text{ such that }\; A^\allyears x^\allyears \leq b^\allyears \;\text{ and }\; \sigma \circ \pi(x^\allyears) = y.
\end{equation}
Given that $\pi$ and $\sigma$ are linear, the above is still a linear program.
Note that $\phi^\allyears$ may not be well-defined for all $y$ if the corresponding linear program has no solutions; this can happen if some total investments $y$ are realisable for each individual weather year, but not realisable by any one complete system design over all 41 years.

Solving \cref{eq:exact-map} amounts to solving an ESOM defined with 41 weather years, which is computationally challenging.
Indeed, one of the motivations for working with the intersection $\Aeps$ is that it can be computed on the basis of model optimisations with single weather years.

Thus, we propose two alternatives to the ``exact'' map $\phi^\allyears$.
We call the exact map and its alternatives \emph{allocations} since they map total investments $y \in \Aeps$ to a spatially resolved allocation of investments $\phi(y) \in \Fepsdash^\allyears$.
The alternative allocations are heuristics in the sense that they map to designs that are not strictly speaking guaranteed to be feasible.
Formally, they map $\Aeps \to \mathbb{R}^{\Ninv}$ but may map some points of $\Aeps$ outside of $\Fepsdash^\allyears$ (whereas $\phi^\allyears(\Aeps) \subseteq \Fepsdash^\allyears$).
For the first alternative, which we call the ``conservative'' allocation, let $i^*$ be the most expensive year, for which $c_{\text{opt}}^{(i^*)} = c^*$.
Then we simply take $\phi^{(i^*)}$ as an alternative to $\phi^\allyears$.
That is, for the conservative allocation we take the solution of the ESOM defined with weather year $i^*$, with the added constraints that $\sigma \circ \pi(x^{(i^*)}) = \ych$.
Computing $\phi^{(i^*)}$ involves only one model optimisation with a single weather year.

For the second alternative, which we call the ``mean'' allocation, we follow the idea of the conservative allocation $\phi^{(i^*)}$, but involve the other weather years more.
Indeed, we define $\phi^{\text{mean}}$ as
\begin{equation}
	\label{eq:mean-allocation}
	\phi^{\text{mean}} = \frac{1}{|\allyears|} \sum_{i \in \allyears} \phi^{(i)}.
\end{equation}
Computing $\phi^{\text{mean}}$ involves $|\allyears|$ model optimisations with single weather years.
In our case, $|\allyears|=41$.
The exact and alternative allocations are summarised in \cref{tab:allocations}.

\begin{table}
	\centering
	\begin{tabular}{llp{6cm}p{4.2cm}}
		\toprule
		Map                  & Short name   & Description                                                                                                                                                 & Computation                                                                  \\ \midrule
		$\phi^\allyears$     & exact        & Fix the investment costs in $\ych$ as additional constraints. Then solve \cref{eq:lin-prog-def} jointly over all years.                                     & Optimisation over 41 years (134~GB RAM, 35 hours)                            \\
		\midrule
		$\phi^{(i^*)}$       & conservative & Fix the investment costs in $\ych$ as additional constraints. Then solve \cref{eq:lin-prog-def} for the year with the highest optimal cost, $c^*$.          & 1 single-year optimisations (ca. 3~GB RAM, 0.25 hours)                       \\
		\midrule
		$\phi^{\text{mean}}$ & mean         & Fix the investment costs in $\ych$ as additional constraints. Then solve \cref{eq:lin-prog-def} for all single years and take the mean over all capacities. & 41 single-year optimisations (parallelisable, each ca. 3~GB RAM, 0.25 hours) \\
		\bottomrule
	\end{tabular}
	\caption{Overview over three different allocations to obtain ``robust'' solutions.}
	\label{tab:allocations}
\end{table}

For comparison, we also define a ``baseline'' point in $\mathbb{R}^{\Ninv}$ which is obtained by scaling up investments uniformly in the optimal solution for the most expensive weather year.
This is explained in more detail in \cref{sec:res-validation}.

\section{Implementation}
\label{sec:implementation}

\subsection{Modelling set-up}
\label{sec:modelling-set-up}

We base our implementation on the PyPSA-Eur 0.4 model \cite{PyPSAEur}, which is itself based on the general PyPSA framework (version 0.18) \cite{brown-horsch-ea-2018}.
While we have modified PyPSA-Eur for our purposes (as described in \cref{sec:data-appendix}), especially in order to support multiple weather years, we have kept the model set-up relatively close to the defaults as described in \cite{PyPSAEur}.
We use the model in a partial greenfield configuration, where existing transmission, nuclear\footnote{Nuclear power in Germany is removed from the model.} and hydro capacities at current (2020) capacities are included in the model from the start, but all other technologies start at zero capacity\footnote{PyPSA-Eur can also include existing biomass capacities in the model. Due to their limited capacities, they do not lead to significant deviations in results compared to the scenario we considered, so we choose to omit biomass for simplicity of the setup.}.
The extendable technologies included in the model are transmission (both AC and relevant DC connections), battery and hydrogen storage, onshore and offshore wind power, solar power and open-cycle gas turbines.
The model is run with a single investment period and perfect foresight.
We limit annual \ch{CO2} emissions to 95\% of 1990 levels\footnote{We present additional results with a 100\% emission reduction in \cref{sec:emission-appendix}.}.
A one-node-per-country\footnote{Except for countries in multiple synchronous zones (Denmark, Spain, Italy, UK), which are represented through two nodes.} spatial resolution and a 3-hourly temporal resolution is chosen.
Note, however, that the spatial and temporal resolution can be increased easily (as in PyPSA-Eur); the resolution is limited in this paper in order to allow extensive validation of our methods.
We model the year 2030 with 41 distinct historical weather years (1980-2020) driving renewable capacity factors and electricity demand based on the default PyPSA-Eur cost assumptions\footnote{\label{footnote:pypsa-eur-costs}  \url{https://github.com/PyPSA/pypsa-eur/blob/v0.4.0/data/costs.csv} (accessed 06/10/2022)} (given in 2013 EUR) for the year 2030.
Thus all weather years are viewed as different potential realisations of 2030, and we can compare investment and operational costs of multi-year optimisations to single-year optimisations by taking annual averages.

Using enough weather data to accurately represent long- and short-term dynamics and extreme events is difficult in ESOMs, considering the resulting model size and increased solving complexity.
However limited the lessons of historical weather data are on future weather \cite{vanderwiel-stoop-ea-2019}, further driven more and more by climate change, the extreme events and variability represented here will still likely offer insights for future designs.
We capture historical climate change implicitly here, whereas incoming trends and changes through climate change are hard to predict and an active field of research in itself \cite{wohland-reyers-ea-2017,schlott-kies-ea-2018,kozarcanin-liu-ea-2019,bloomfield-brayshaw-ea-2021}.

\subsection{Data}
\label{sec:data}

The data we use as input for our model lean heavily on the data and sources used in PyPSA-Eur.
For instance, we use cost data based on the default costs considered in PyPSA-Eur, which are estimates for 2030\textsuperscript{\ref{footnote:pypsa-eur-costs}} (in 2013 EUR).
However, time series input data have been extended to 41 weather years.

We use reanalysis data to generate capacity factors for the renewable energy sources; our source is the hourly ERA5 dataset \cite{era5-data} for the time period from 1980 up to and including 2020.
The open-source tool Atlite \cite{hofmann-hampp-ea-2021} translates weather data to hourly capacity factors for solar PV, on- and offshore wind.
The inflow profiles for hydropower are generated similarly, however they are corrected to fit production data given by \cite{eia-hydro-data}.
We describe in detail how we process the data in \cref{sec:data-appendix}.

The load input data are generated using a regression model on 1980--2020 ERA5 temperature data.
The regression is based on hourly country-level ENTSO-E load data from 2010 to 2014 \cite{entso-e-2022}, as well as the temperature data.
We conduct a two-staged regression with a similar approach to that used in \cite{bloomfield-brayshaw-ea-2020}.
This allows us to generate 41 years of country-level temperatur-dependent load data matching the weather data we use.
Finally, we scale the demand by a factor of 1.13 according to load projections for 2030 by the European Commission\footnote{\url{https://ec.europa.eu/clima/document/download/ec1acac9-10fe-4eeb-915f-cad388990e0f_en}, Fig. 44 (accessed 23/06/2022)}.
More details on the generated load data can be found in \cref{sec:data-appendix}.

\subsection{Modeller's decisions}
\label{sec:modellers-decisions}

In this section, we discuss various details and choices regarding the implementation of the methods described in \cref{sec:methods-formal-defs}.
For the exact code, including configuration options and installation and running instructions, we refer to the GitHub repository\footnote{\url{https://github.com/aleks-g/intersecting-near-opt-spaces/tree/v1.0.1}}.

One of the first decisions we have to make is choosing a suitable slack level $\varepsilon$.
For small $\varepsilon$ the intersection $\Aeps$ may be empty; a priori it is not clear how large $\varepsilon$ has to be for $\Aeps$ to be nonempty.
In our case, we choose $\varepsilon = 5\%$ on top of the most expensive weather year, but found that $\Aeps$ is even nonempty with $\varepsilon = 2.5\%$; this may change with a different modelling set-up.
For comparison, Trutnevyte found in \cite{trutnevyte-2016} that the transition in the UK energy system from 1990 to 2014 deviated between 9 and 23\% from the cost optimum.

Given that the dimension reduction map $\sigma$ is our main tool in working with near-optimal spaces (see \cref{sec:dim-reduction}), we need to define it carefully.
Both the number of dimensions $k$ that $\sigma$ maps to and which investment decision variables are mapped to each dimension must be considered.
We have investigated the convergence of \cref{alg:near-opt-approx} with $k = 2, 3, \dots, 7$ (\cref{sec:direction-generation}) and we find that it is tractable to work with this number of dimensions.
The decision variables that are mapped to each dimension (the sets $T_1, \dots, T_k$ in the notation of \cref{sec:dim-reduction}, here abbreviated to ``dimensions'') should be chosen meaningfully.
On one hand, including as a dimension a technology which is not utilised in cost-optimal solutions may be detrimental, as much computational effort will be expended on potentially irrelevant solutions including this technology.
Moreover, the Chebyshev centre of the resulting space must include at least a Chebyshev radius worth of that technology, which may be sub-optimal.
On the other hand, \emph{not} including as a dimension a technology which plays a significant role in any near-optimal solution can limit the usefulness of the results.
In our case, we choose to map investment decision variables for transmission expansion, solar, onshore wind, offshore wind and gas to 5 respective dimensions.

Lastly, we choose the coefficients $c_j$ in the definition of $\sigma$ to be the investment cost associated with the investment decision variable $x_j$; this has the advantage of mapping to the single unit of EUR in every dimension, making different dimensions easily comparable.
In contrast, for example, comparing renewable capacity expansion in MW and storage capacity expansion in MWh directly is more difficult and not as useful for the purpose of working with a Chebyshev ball.
For some applications one might consider scaling the weights $c_j$ for certain dimensions, for example, when robustness to changes in investment in one technology are more important than for other technologies.

A key driver of the computational demands of our approach is the desired quality of approximation of near-optimal spaces.
How many iterations of \cref{alg:near-opt-approx} are needed, depends on the number of dimensions $k$, the mode of finding new directions to explore (see \cref{fig:direction-performance}) and the intended use-case.
A fixed number of iterations can be chosen, or the algorithm can be ended once some convergence criteria is satisfied.
Either way, we refer to \cref{sec:direction-generation} for a discussion of the various trade-offs involved.

The other factor influencing the computational effort is how much time and computing power every single optimisation takes.
This is typically driven by the number of technologies, inter-temporal relations between different variables (e.g. through storage), temporal resolution \cite{hoffmann-kotzur-ea-2020} and spatial resolution \cite{trondle-lilliestam-ea-2020, frysztacki-horsch-ea-2021}.
Deciding on the model complexity must be done in light of the research question at hand.
For our application with weather years, we do assume that each individual model is defined over a time period of one calendar year.

\section{Results}
\label{sec:results}
We first present the main results pertaining to system design with 41 weather years, and the use of near-optimal spaces in this context (\cref{sec:res-weather-years-feasibility}).
This is followed by a more detailed description of the characteristics of our proposed robust solutions (\cref{sec:res-robust-design}).
The subsequent subsection focuses on validation results and a comparison of the different robust design allocations (\cref{sec:res-validation}).
Finally, we touch on computational results and parallelisability (\cref{sec:res-computation}).

\subsection{Weather years and intersection}
\label{sec:res-weather-years-feasibility}

First of all, optimising our model with each of the 41 considered weather years individually shows large discrepancies in the respective optimal solutions, re-affirming the importance of considering a large set of different weather years.
Optimal total system costs range from 121 billion EUR for 2020 to 152 billion EUR in 1985\footnote{These and all the following total system costs are annualised and include investment in new capacities as well as variable costs (in 2013 EUR), but not existing capacities. See \cref{sec:modelling-set-up} for details on which existing technologies are included in the model.}.
The composition of investment by technology also differs significantly between weather years, with especially the onshore- and offshore wind investment varying by up to around 20 and 25 billion EUR between years (corresponding to variations up to 243 GW for onshore wind and 82 GW for offshore wind), respectively.
The inter-year variability in optimal investments is illustrated in \cref{fig:over-under-investment}, where the investments are compared to the robust solution $\xrob$.

Meanwhile, for this study we also conduct the first capacity expansion optimisations of a spatially resolved model for the European power system with 41 weather years directly (one node per country, 3-hourly resolution).
These optimisations took in the order of 1-2 days (using two threads) and up to 134~GB of memory.
The optimal annualised total system cost for the 41-year model is 137 billion EUR; only slightly higher than the average of the total system costs of optimisations with single weather years at 134 billion EUR.
See also \cref{fig:opt-costs} for a break-down into investment (for extendable technologies)- and variable costs and comparison with the optimisations with single weather years.
Apart from providing a basis of comparison, optimising the model with 41 weather years also allows us to compute the exact robust map $\phi^\allyears$ in order to validate the alternative mean and conservative allocations (\cref{sec:res-validation}).

We implement (as in \cref{sec:implementation}) the methodology of intersecting near-optimal spaces (laid out in \cref{sec:methods-formal-defs}) using PyPSA-Eur.
With a slack level of $\varepsilon = 0.05$ as in \cref{eq:uniform-slack}, we obtain a nonempty intersection $\Aeps$ and a robust allocation $\xrob$ that is fully feasible using all the weather data of 41 years.
Thus we show that there are robust solutions with less than 5\% additional costs (on top of the most expensive year).
We even find robust solutions which are less than 5\% more expensive than a system optimised with the entire period of 41 years (see \cref{fig:opt-costs}).
Furthermore, the space $\Aeps$ offers significant flexibility for policymakers beyond our point of reference $\xrob$, and beyond what flexibility is shown by previous MGA approaches.

\begin{figure}
	\centering
	\includegraphics{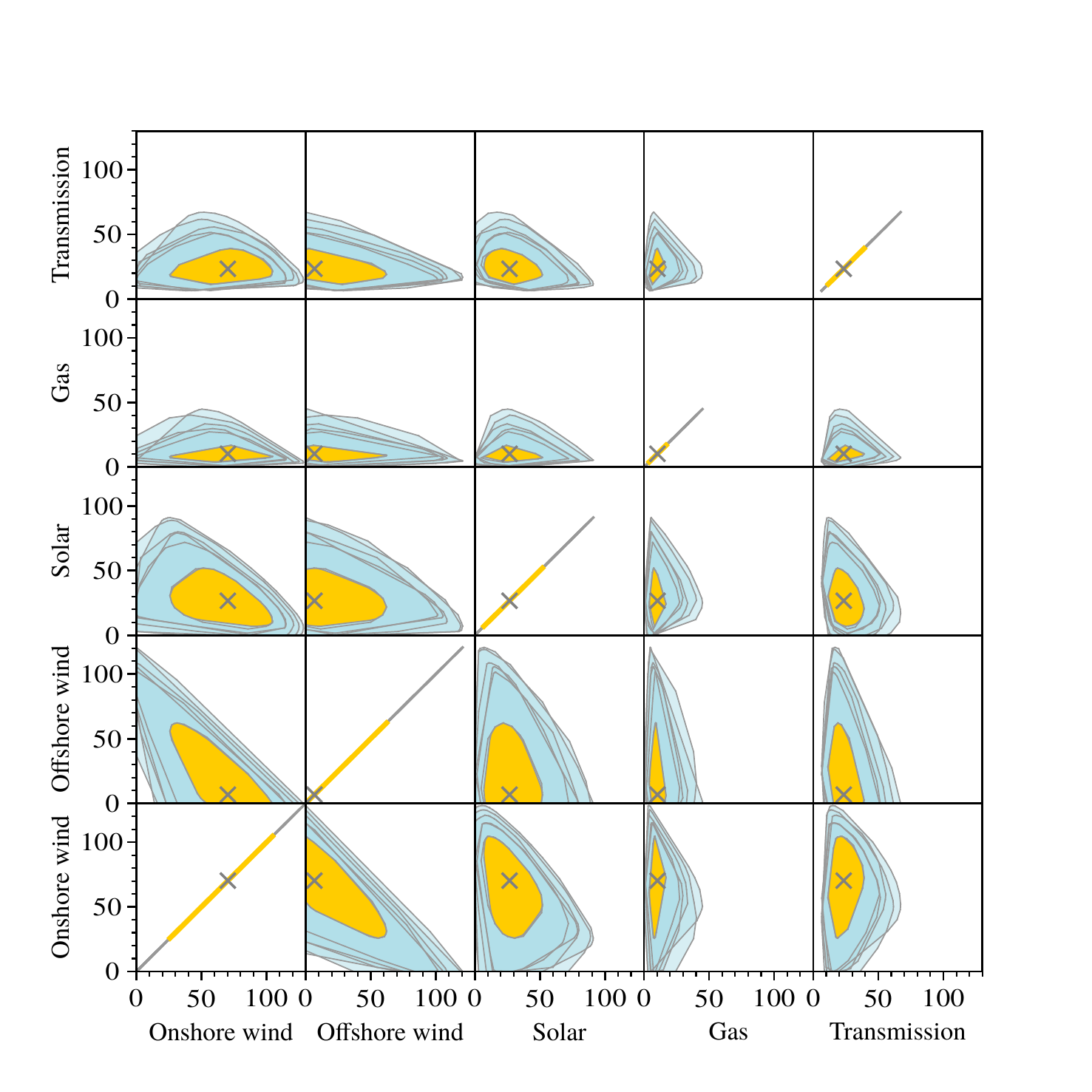}
	\caption{Projections of the near-optimal spaces for different weather years and their intersection. All values are annualised total investment costs per technology. For illustrative purposes, we only plot the near-optimal spaces for 6 out of 41 weather years (in different hues of blue). The intersection of all 41 near-optimal spaces is filled in yellow and the Chebyshev centre is marked with a cross.}
	\label{fig:near-opt-grid}
\end{figure}

\Cref{fig:near-opt-grid} shows projections of a selection of near-optimal spaces $\Aeps^{(i)}$ for weather years $i \in \{1985, 1989, 1996$, $2006$, $2014, 2020\}$ in addition to the intersection $\Aeps = \bigcap_{i \in \allyears} \Aeps^{(i)}$ over all weather years.
While the spaces are 5-dimensional (with the dimensions being investment in onshore wind, offshore wind, solar, gas and transmission expansion), they have been projected down to all possible pairs of technologies considered.
The Chebyshev centre, marked by a cross, is located within the intersection which consists of the robust solutions.
The figure reveals that there is significant flexibility in these dimensions; while a certain amount of investment in renewables is needed, the investment can be shifted between different technologies while staying feasible and near-optimal.
Note also that the near-optimal spaces for different years resemble each other in shape and location in space, but mainly differ in size.
This indicates that the effect of ``difficult'' weather years on modelling is mainly that they restrict the size of the feasible design space.

We find that the Chebyshev radius of $\Aeps$ is $3.43$ billion EUR, coming near to the theoretical maximum possible radius of $3.80$ billion EUR given by the chosen slack level.
Indeed, note that the distance between cost-optimal solutions and the near-optimal cost constraint is $c_{\text{opt}}^* \cdot \varepsilon$, meaning that any near-optimal space can have a Chebyshev radius of at most $c_{\text{opt}}^* \cdot \varepsilon / 2 \approx 3.80$ billion EUR.
The result means that the total investments in technologies that make up the dimensions of $\Aeps$ can change by up to 3.43 billion EUR (corresponding to 2.35\% of the total cost of the robust system) in any direction, starting at the robust point $\ych$.
The resulting (potentially reduced) total investments can still result in a feasible design for every weather year.

Given that we work with an absolute objective bound (1.05 times the cost of the most expensive year) for all near-optimal spaces, we find that the smallest and largest near-optimal spaces differ in volume by a factor of about 79.
This means that some weather years by themselves allow for many more different near-optimal feasible solutions than others.
Put differently, some weather years restrict the system design much more than others.
The smallest and largest near-optimal spaces come from the weather years 1985 and 2020 respectively, and the difference in volume corresponds to an average scaling factor of $79^{1/5} \approx 2.4$ in every dimension.
Meanwhile, the intersection of the near-optimal spaces has a volume that is 79\% and 36\% of the volumes of the near-optimal space for years 1985 and 1987, respectively, and is between 1\% and 10\% of the volume of all other near-optimal spaces.
This means that except for 1985 and 1987, 10\% or less of near-optimal solutions for any particular weather year are feasible for all other weather years under consideration.

In fact, we find that the optimal solutions with the years 1985 and 1987 actually have total investments that lie within the intersection $\Aeps$.
When operated over the entire weather year dataset, these designs are practically feasible, with negligible load shedding.
These results can inform the choice of weather year to model with --- if only a single (or few) years can be chosen.

We note, however, that while there are weather years to which our model does not ``over-fit'' in a single-year optimisation (1985 being such a year resulting in a generally applicable design), this result could be particular to our modelling set-up.

\subsection{Robust design characteristics}
\label{sec:res-robust-design}

The robust solutions we compute, coming from the Chebyshev centre of the intersection $\Aeps$, have several interesting properties.
First, we compare the investment composition of the robust point $\ych$ to that of optimisations using just a single year of weather data.
We then analyse total investment and operational costs.
Finally, we present results related to the \ch{CO2} limit.

\begin{figure}
	\centering
	\includegraphics{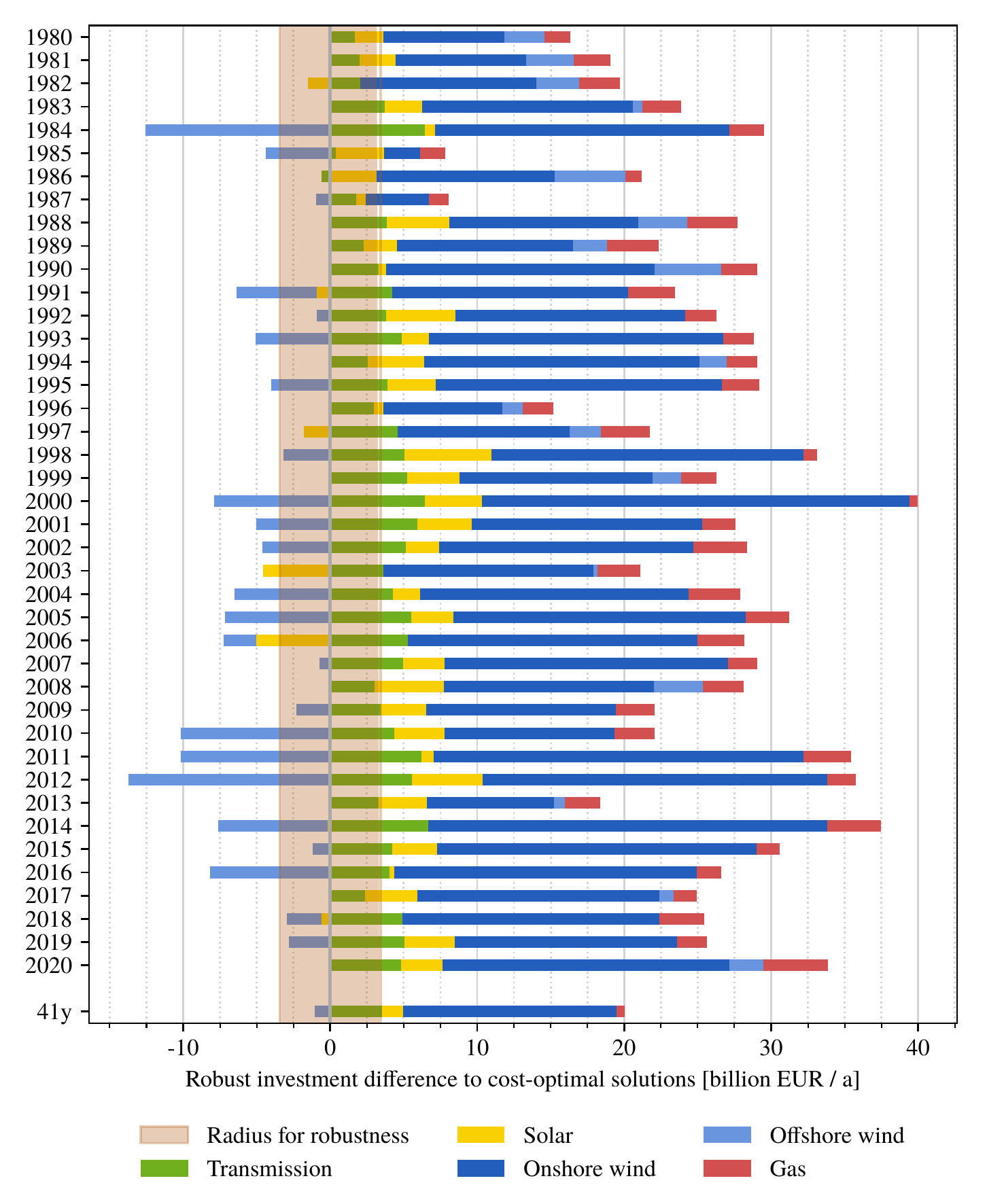}
	\caption{Comparison of total investments for selected technologies in the optimal solutions for each weather year and the optimal solution with all weather years (``41y'') to the robust point $\ych$. Positive values mean greater investment in the given technology by the robust solution.}
	\label{fig:over-under-investment}
\end{figure}

\Cref{fig:over-under-investment} shows the differences in investment per technology between the robust point and each of the optimal designs for individual years.
The robust point is characterised by more investment in onshore wind, solar, transmission, and gas (sorted by decreasing additional investments).
Meanwhile, most optima from single weather years over-invest significantly in offshore wind power compared to the robust allocation due to higher relative costs.
In this particular set-up where gas can smooth the electricity production, the cost benefits of additional onshore wind capacities outweigh the potential of offshore wind power in the most favourable years.
We conclude that here onshore wind power contributes more than other technologies to robustness, followed by solar power and transmission capacity\footnote{The conclusions can change with different assumptions --- see \cref{sec:emission-appendix}.}.
This holds as well when one compares the investments in the robust (``exact'') allocation to the optimal solution using all weather years (see the ``41y'' row in \cref{fig:over-under-investment}).

\begin{figure}[tbh]
	\centering
	\includegraphics{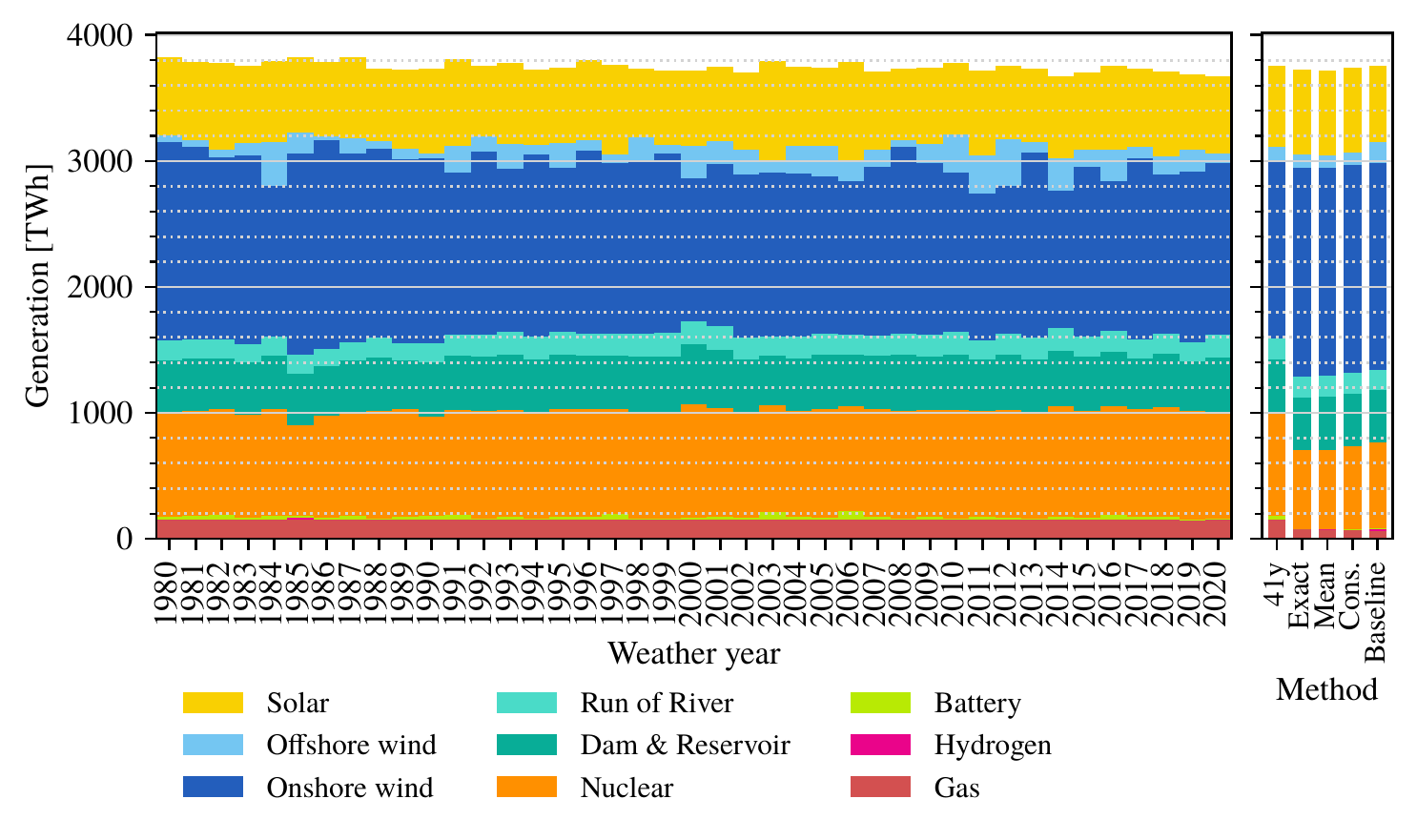}
	\caption{Annual average generation mixes for the optimal solutions for each single weather year, the optimal solution with all weather years (``41y'') and the robust allocations as well as the baseline design. The baseline design is used for comparison in validation (\cref{sec:res-validation}).}
	\label{fig:generation}
\end{figure}

\Cref{fig:generation} shows the annual (for the optimisations with more weather years, average) total electricity generation per technology.
This shows the increased importance of onshore wind and solar for the robust system.
Meanwhile, the figure also shows that although solutions for individual years typically under-invest in gas turbines relative to the robust ones, the robust solutions actually generate less power with gas (and nuclear) in total.
This reflects the fact that while additional gas capacity is needed for robustness, the additional investment in renewables leads to a reduced dependence on gas for ``day-to-day'' operations.

\begin{figure}[tbh]
	\centering
	\includegraphics{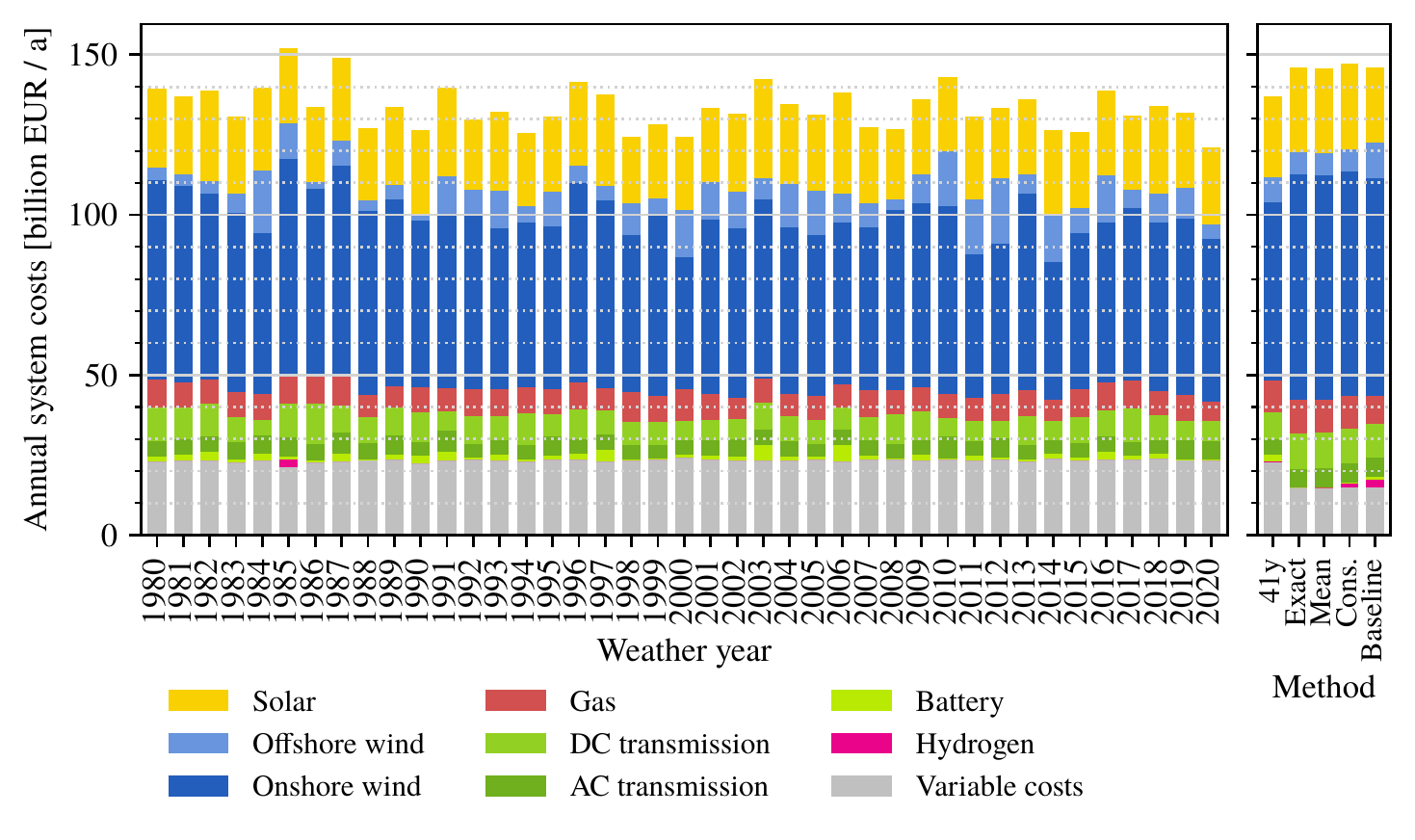}
	\caption{Comparison of cost-minimal designs based on optimisations over individual years to the annualised costs of robust designs (exact, mean and conservative, respectively), an optimisation with all 41 years (``41y'') and the baseline design.}
	\label{fig:opt-costs}
\end{figure}

On that topic, \cref{fig:opt-costs} shows the variable and total system costs of systems optimised with single weather years, as well as the robust system and a system optimised with the full dataset of 41 weather years.
It illustrates that the robust design has a higher total investment cost than designs for individual years, while the (average) operating costs are lower due to the reduced use of gas and nuclear as mentioned above.
Recall that while we set the slack $\varepsilon$ to 5\%, we see that the investment costs in fact lie only 0.4\% above what would have to be invested based on the most expensive year.
With the (average) variable costs of the robust system being lower due to a strengthening of renewables, the total (annualised) system cost of the robust solution is about 146 billion EUR and actually lower than the system costs for some optimal solutions with single weather years.
This is because the total system cost for the robust solution is averaged over all 41 weather years, with some being more expensive than others.

We also see that robust system designs emit less \ch{CO2} in our tests compared to the single-year optimisations, and use 48\% of the given \ch{CO2} limit\footnote{The optimal solution with all 41 weather years uses the whole \ch{CO2} limit as well.}.
This is again because robust designs direct more of the total system cost into capital investment of renewables and less into variable costs including gas, the only source of emissions in our model.
However, we should note that when we operate the design obtained from e.g.\ a system optimisation with the single weather year 1985 over the entire time period (over which it is practically feasible), it also does not use up the whole \ch{CO2} limit.
Although capacity expansion optimisations with single weather years always use up the \ch{CO2} limit, the designs which are adapted to difficult years such as 1985 have enough renewable capacities that they do not use up the \ch{CO2} limit in a typical year.

\begin{figure}[tbh]
	\centering
	\includegraphics{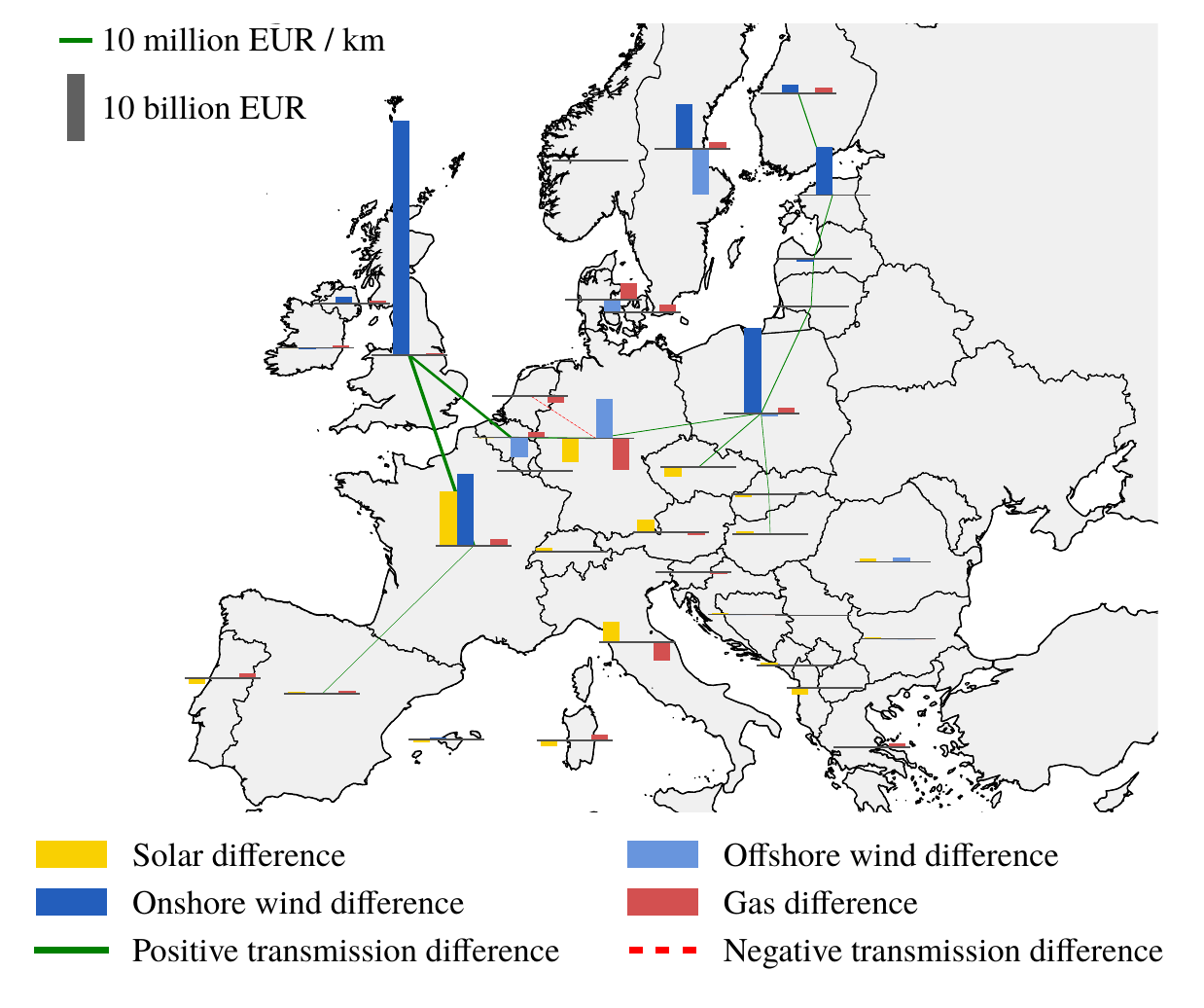}
	\caption{Difference in investment between the exact robust allocation and the optimum solution with 41 weather years. Bars above the baselines mean that there is more investment for that technology in the robust solution. For transmission there are only positive differences (negative differences are too small to appear on this scale), and AC and DC transmission have been combined. For both bars and transmission lines, the same area means the same investment difference in EUR. All costs are annualised.}
	\label{fig:map-difference}
\end{figure}

Finally, \cref{fig:map-difference} shows the geographical differences in investment between the exact robust allocation and the cost-optimal solution with 41 weather years.
We see that the majority of additional onshore wind power in the robust system is allocated to the UK, followed by Poland, France, Estonia and Sweden.
The transmission capacity to the UK also receives more investment in the robust allocation.
Meanwhile, although the system invests less in offshore wind in total, more offshore wind investment is shifted to Denmark.

\subsection{Validation of robust allocations}
\label{sec:res-validation}

Recall from \cref{sec:methods-formal-defs} that we can use a (large) number of optimisations with single weather years to find a robust point of \emph{total} capacity allocations $\ych \in \Aeps$ --- the Chebyshev centre in the intersection of near-optimal spaces of model instances using different weather years.
However, in order to find \emph{specific} (per-node) investment allocations fitting the given robust totals and being feasible for every weather year, we have to map $\ych$ back to $\Fepsdash^\allyears$.
The exact allocation $\phi^\allyears$ does this by solving the original ESOM (with the additional constraint that $\sigma \circ \pi(x) = \ych$) with 41 weather years.
\Cref{sec:going-back} gives two alternatives for allocating the robust capacities to individual nodes: the ``mean'' and ``conservative'' allocations, both based only on a number of optimisations with one single weather year at a time (see \cref{tab:allocations}).
Recall that both the mean and conservative allocations have the same coordinates in $\Aeps$ as the exact allocation, namely $\ych$, but are not guaranteed to be feasible for all 41 weather years.

In this section, we investigate the quality of these heuristic allocations.
The basis for comparison is two-fold: the exact robust allocation on the one hand, and a ``naïve'' baseline allocation based solely on the most expensive weather year on the other hand.
Specifically, in order to make the designs more comparable, the baseline allocation is obtained by taking the allocation from the single most expensive weather year (1985 in our case) and uniformly scaling all expandable capacities up by such a factor that the total capital cost of the whole network equals that of the exact robust design.
This way, the only difference between the exact, mean, conservative, and baseline designs is how investment is allocated (by technology and spatially), not the total investment volume.

The validation of the robust allocations consists of a ``stress test'' where we operate the systems over the entire dataset of weather years (i.e. only optimise the dispatch and not the investments).
This dispatch optimisation includes the 95\% \ch{CO2} reduction constraint as in the capacity optimisation.
Similarly to how we ensure that the total investment cost is equal between the four different allocations, we also ensure that all designs keep to the same operational budget; namely the total operational costs of the exact robust design.
By allowing costly load shedding (with variable costs of 7 300 EUR/MWh, as in \cite{price-zeyringer-2022a}), we make sure that the models are solvable.
The quality of the allocations is then measured in the amount of load shedding; 0 load shedding means complete feasibility.
In the presence of a hard operational budget constraint, it should be understood that total unmet demand combines actual unmet demand due to technical infeasibility and over-budget operations.
Without the operations budget constraint, a number of designs --- including all robust allocations and also some optimal solutions for single years --- are practically feasible with all weather years (see \cref{sec:res-weather-years-feasibility}), but simply have higher operational costs.
In general, however, optimal solutions from single years fail to serve the entire time period reliably.

\Cref{tab:validation} shows that the exact allocation has 0 load shedding as expected, and is followed by the conservative allocation in quality.
The mean and baseline designs perform slightly worse, but still only up to about 0.1\% of the total load is shed or produced over budget.
This shows that the conservative and mean heuristics produce results that are practically feasible.

\begin{table}[th]
	\centering
	\begin{tabular}{lrr}
		\toprule
		Allocation method & Total load shedding [TWh] & Relative load shedding [\%] \\ \midrule
		Exact             & 0.0                       & 0.000                       \\
		Conservative      & 46.7                      & 0.032                       \\
		Mean              & 119.0                     & 0.081                       \\
		Baseline          & 132.3                     & 0.090                       \\
		\bottomrule
	\end{tabular}
	\caption{Performance of the different allocation methods and the baseline for comparison. Note that total load shedding (over 41 years) here includes technical infeasibility and over-budget operations as explained earlier. The relative load shedding is with respect to the total load of all 41 weather years.}
	\label{tab:validation}
\end{table}

\subsection{Performance and computational effort}
\label{sec:res-computation}

Our work demonstrates methods for working with decades of weather and demand data in a parallelised setting.
Namely, the mean and conservative allocations produce robust system designs with decades of weather data while only requiring model optimisations with single weather years which can be parallelised effectively.
This is done by using the centre point of the intersection of near-optimal spaces for different weather years to obtain a robust set of total capacities per technology.
From these we then compute specific per-node capacity allocations adhering to the robust totals for each individual weather year, and average the resulting capacities.

This is particularly relevant because solving ESOMs is usually memory-constrained: a single large model may take more memory to solve than is available on common systems.
While effort is being spent on splitting up and parallelising the solving process \cite{rehfeldt-hobbie-ea-2022}, this is not yet practical.
With our methods we avoid solving models defined over a period of decades, and instead obtain results based on many smaller runs --- a process which is easily parallelised.
The workflow is now constrained by the number of available processors, which is better suited to current computational developments which are in the direction of more, not faster processing cores.

While computational time often varies significantly based on a multitude of factors (exact modelling set-up, model size, nature of objective and constraints, numerical issues among others) and is difficult to predict, the memory requirements are easier to derive directly from model size (number of technologies, temporal and spatial resolution).
With our modelling set-up (see \cref{sec:modelling-set-up}), a capacity expansion optimisation with a single weather year takes 3.5--4~GB of memory using the commercial optimisation software Gurobi \cite{gurobi}.
Thus, our methods put optimal energy system design with decades of weather data within reach of typical desktop computers for the first time.
For comparison, optimising with 41 years of weather data in one model takes approximately 41 times as much memory, around 134~GB.
For models with a higher spatial and temporal resolution, our methods may currently be the only viable method of designing a system with decades of weather data; we estimate that optimising a typical PyPSA-Eur model with an hourly resolution and 181 nodes (as suggested in \cite{frysztacki-horsch-ea-2021}) with 41 weather years could take around 1.7TB of memory, while this can be split into single-year optimisations taking around 45~GB of memory each using our methods.

With our modelling set-up, a model optimisation with one weather year takes 23 minutes of CPU time\footnote{\label{cpu-wall-time}We make a distinction between the \emph{CPU time} of a process, which is the sum of time spent on the process by all CPUs, and the \emph{wall time} of a process, which is the elapsed real time between process start and finish. Thus CPU time may be larger than wall time for a parallelised process.} on average using Gurobi with 2 threads.
Spending 10 initial optimisations along unit vectors in positive and negative direction, 150 for the algorithm to approximate the near-optimal space for each year and finally 1 optimisation per year to compute per-node allocation of robust allocations, the mean allocation takes approximately $41 \cdot (10 + 150 + 1) \cdot 23 = 151823 \text{ min.} \approx 2530 \text{ hours}$ of CPU time to compute.
On a cloud computing platform with 64 CPU cores (and around $4 \cdot 64 = 256$~GB of memory), this equates to around 40 hours of wall time\cref{cpu-wall-time}.
This is actually comparable to the wall time it typically takes to solve a single capacity expansion optimisation with 41 weather years of the same model on the same platform, which is in the range of 1-2 days.

We also investigate different methods for approximating near-optimal spaces of ESOMs, focusing on how optimisation directions are chosen, and how many optimisations are needed for a good approximation (for a different number of dimensions or projection variables).
The results presented in \cref{sec:direction-generation} support the configurations we choose as default options.

\section{Discussion}
\label{sec:discussion}

The urgent transition to energy systems based on intermittent renewable generation comes at a time with increasing computational power and availability of extensive climate data.
However, while using ever larger models is helping us to understand the detailed functioning of future energy systems, it does not necessarily improve the understanding of uncertainties and resilience.
In this paper, we propose a framework for producing energy system designs that are robust to uncertainty, taking advantage of the geometry of near-optimal spaces.
We apply the framework to a study involving 41 weather years; the importance of using as much weather data as possible has been shown in the literature and is confirmed by our findings.
Our framework can help policymakers to investigate different alternatives (similarly to the predecessors of this methodology, MGA and MAA) and overcome the fallacy of single solutions that can become infeasible through marginal perturbations.

Our choice of using the Chebyshev centre as a robust point is motivated by it being maximally tolerant to changes in investment in any direction.
This is one of the distinctions to robust optimisation, where the optimum solution (with respect to the worst-case scenario) lies at the boundary of the intersection of feasible spaces, thus becoming infeasible with just marginal perturbations.
At the same time, we see that even a naïve over-investment starting at the cost-optimal solution with a difficult weather year results in a system design which is practically feasible for all considered weather years.
The choice of the Chebyshev centre can be seen as an attempt to find the most advantageous over-investment for robustness.
In other words, while all additional investment contributes to robustness to some extent, our methods are aimed at finding the most efficient additional investments for robustness.

The first runs of a spatially resolved European power system optimisation model with four decades of weather data build on advances in open data accessibility and computational resources.
Our results show that robust solutions stand out from standard cost-optimal solutions by a relatively higher investment in onshore wind and solar power.
Furthermore, the increased investment in renewables reduces the usage of gas and nuclear, which lowers the \ch{CO2} emissions of the system by more than 50\% (in comparison to optimising the system with 41 weather years).
By adjusting the slack level $\varepsilon$, our methods allow for a trade-off between robustness and additional costs.

Beyond looking at how weather years influence total system cost or investment composition, we see the effects of the choice of the weather year on the entire near-optimal space.
Whereas the shapes and locations are fairly uniform, we see that near-optimal spaces for different weather years vary significantly in size.
We still find a significant amount of flexibility within their intersection, mainly limited by a small number of most difficult weather years.

In previous studies utilising MGA techniques, dimension reduction has been used somewhat implicitly in order to make a systematic exploration of reduced near-optimal spaces feasible.
Mapping down to a lower-dimensional space was seen as a way of summarising a few key properties of model solutions.
However, the \emph{geometry} of reduced near-optimal spaces has been largely unexploited, a gap we aimed to fill with this paper.
Moreover, we investigate and propose practical methods for going back to spatially explicit investment decisions from points in the reduced near-optimal space.
Mapping the Chebyshev centre of an intersection of near-optimal spaces back to a robust system design is just one application of this idea.

At the same time, approximating and intersecting near-optimal spaces as a way of designing systems which are feasible for many weather years is a reasonable alternative to optimising with many weather years directly from a computational perspective.
Having split the problem into many smaller optimisations with single weather years, our methods are easily parallelised.
This represents a new approach to optimisation problems which previously have been considered intractable or at least very impractical.
In particular, an application of our methods to larger models (including sectoral coupling, higher spatial or temporal resolution) would be interesting and computationally feasible.
While we expect the heuristic allocations to work at higher spatial resolutions (see also \cite{frysztacki-hagenmeyer-ea-2022} for related work), their quality should be investigated under such conditions.

We highlight that we study just one out of many definitions of robustness \cite{moret-bierlaire-ea-2016, maggioni-potra-ea-2017}; note also that our techniques do not directly follow the concept of \emph{robust optimisation}.
In all generality, robustness is a relative concept and is directed towards some uncertainty.
Although these uncertainties can be well understood, some of them may be hard to quantify (e.g. political, societal changes), other could be epistemic (e.g. extreme weather events, changes in costs, or misspecifications of the model), or even aleatoric (as the actual future weather conditions).
For instance, Stirling \cite{stirling-2010} locates robustness in the overlap of problematic levels of knowledge about possibilities and probabilities.
All in all, our approach contributes to a wide academic debate about how to deal with uncertainty and robustness in energy systems modelling.

\section{Conclusion}
\label{sec:conclusion}
In this article we find that studying the near-optimal feasible space is a helpful tool to achieve more robust solutions against uncertainties for energy systems.
When investigating weather variability, this enables us to quantify the variations and to find alternative designs that allow some flexibility for policymakers.

Since we utilised historical climate reanalysis data, we have not incorporated the consequences of climate change that an energy system in transition will face, nor the possibility of unseen extreme events.
It would be interesting to understand which developments will limit flexibility and what events are defining for the design of a future (climate-)robust energy system.
Moreover, we see that a small number of weather years including 1985 and 1987 constrain our system design the most, having relatively small near-optimal feasible spaces.
These results call for a deeper understanding of which meteorological properties of these weather years, including extreme and compound events, are determining for energy system design.

While we have applied a notion of robustness in the reduced near-optimal space $\Aeps$, we have not considered robustness at lower levels of the network.
So although the point $\ych \in \Aeps$ is robust to a shift in investments of 3.4 billion EUR, the system $\xrob \in \Fepsdash$ is not robust to such a shift at any one particular node in the network.
There may also be many system designs with very different per-node capacities mapping to the same point $\ych$.
Another interplay with the spatial dimension is whether different regions contribute differently to the robustness of the whole system, as hinted at by \cref{fig:map-difference}.
Thus, spatial aspects of robustness form an interesting avenue for future research.

More generally speaking, our methods leave a lot of room for different choices of dimension reduction: choosing which variables to aggregate and how.
The choice leading to the most suitable reduced near-optimal space should be considered application-specific, depending on the technologies of interest, their relative importance and use-case for the reduced near-optimal space.
Mapping to a space of investment costs as we did is a neutral choice, but which reduction is the most efficient for more specific purposes is still an open question.

The results we present here may depend on our particular modelling set-up.
To which extent specific technologies contribute to robustness and flexibility may change if other technologies are included (e.g. through sector coupling), additional restrictions are introduced (e.g. on transmission) or the spatial and temporal resolutions of the model are improved.
And whereas our model has a single investment period, it would be interesting to apply our methods with a model considering transition pathways through multiple investment periods.
As our implementation is open-source and customisable, it should be adaptable to this setting.

Last but not least, our methodology can also be used to investigate other uncertainties besides weather variability.
As near-optimal spaces strongly depend on cost assumptions, future applications of the present framework can contribute to an improved understanding of robustness in the face of uncertain costs.

\section*{Acknowledgements}

We would like to thank the anonymous reviewers for their helpful comments and suggestions.

The computations were partly performed on the Norwegian Research and Education Cloud (NREC), using resources provided by the University of Bergen and the University of Oslo. \url{http://www.nrec.no/}.

A.G., F.E.B. and M.Z. acknowledge funding by UiO:Energy (SPATUS).

\section*{Code \& data availability}
The code for our approach and its documentation can be found at \url{https://github.com/aleks-g/intersecting-near-opt-spaces/tree/v1.0.1}, and is made available under the GPL 3.0 license.
Some of the data included in the above repository (load and hydropower capacity factors) were generated using the workflow available at \url{https://github.com/aleks-g/multidecade-data/tree/v1.0} (code licensed GPL 3.0, data licensed CC-BY-4.0).

The reanalysis data we used for this study (for renewable capacity factors and temperature-dependent load data) were downloaded from the Copernicus Climate Change Service (C3S) Climate Data Store~\cite{era5-data} using the Atlite software~\cite{hofmann-hampp-ea-2021} as documented in the main repository.
We have additionally made these weather data more easily available at \url{https://doi.org/10.11582/2022.00034}\footnote{See instructions at \url{https://github.com/aleks-g/intersecting-near-opt-spaces/tree/v1.0.1}}; they are shared under the ``License to use Copernicus Products''\footnote{\url{https://cds.climate.copernicus.eu/api/v2/terms/static/licence-to-use-copernicus-products.pdf} (accessed 23/06/2022)}, which is comparable to the CC-BY license.
Neither the European Commission nor the European Centre for Medium-Range Weather Forecasts is responsible for any use that may be made of the Copernicus information or data it contains.

Given the substantial size of the weather data needed for this study, we generated a number of PyPSA networks which can be used to reproduce our results without having to download and process the ERA5 data.
They are available at \url{https://doi.org/10.5281/zenodo.6683829} under the CC-BY-4.0 license.

All other data used in our model is directly inherited from PyPSA-Eur and also openly available as described in~\cite{PyPSAEur}.

\section*{Author contributions}
\textbf{Aleksander Grochowicz}: Conceptualisation, Data curation, Formal analysis, Investigation, Methodology, Software, Validation, Visualisation, Original draft \& editing;
\textbf{Koen van Greevenbroek}: Conceptualisation, Data curation, Formal analysis, Investigation, Methodology, Software, Validation, Visualisation, Original draft \& editing;
\textbf{Fred Espen Benth}: Conceptualisation, Methodology, Supervision, Review \& editing;
\textbf{Marianne Zeyringer}: Conceptualisation, Methodology, Supervision, Review \& editing

\section*{Declaration of competing interest}
None.

\bibliographystyle{elsearticle-num}
\bibliography{bibtex-bibliography}


\appendix

\section{Direction generation}
\label{sec:direction-generation}
\setcounter{figure}{0}
\renewcommand{\thefigure}{A.\arabic{figure}}
We give a more detailed overview on how to explore the near-optimal space $\Aeps$ (here in the original sense as a reduced near-optimal space and not the intersection), using the notation from \cref{sec:methods-formal-defs}.
Recall that we can compute an approximation of $\Aeps$ by finding a number of its vertices (or rather, extreme points), and each such point is obtained by solving the linear program in \cref{eq:new-objective} with some different objective (\emph{direction}) $d$.
As described in \cref{alg:near-opt-approx}, we first optimise over $\Aeps \subseteq \mathbb{R}^k$ in each of the cardinal directions (positive and negative) in order to obtain a first full $k$-dimensional approximation.
Each of these optimisations can be performed in parallel.
Thereafter, in which direction $d$ we choose to optimise over $\Aeps$ has a significant effect on how well $\Aeps$ can be approximated in a limited number of optimisations.

We investigate and compare three methods that generate directions in which to optimise over $\Aeps$.
The first method is the simplest and consists of choosing directions uniformly at random.
For the remaining two methods, the idea is to compute the convex hull of the points obtained so far after every optimisation, say $H_i$, and use the geometric properties of $H_i$ in some way to generate the ``next'' direction $d_{i+1}$.
The convex hull is computed using the program qhull \cite{barber-dobkin-ea-1996}.
The three methods are as follows:
\begin{enumerate}
	\item ``random-uniform'': Choose a random vector from the uniform distribution on the sphere in $\mathbb{R}^k$.
	\item ``facets'': Choose the normal vector to the facet of $H_i$ with the largest volume.
	\item ``maximal-centre-then-facets'': Compute the Chebyshev centre $\ych$ of $H_i$ and the ball of maximal radius around it by solving \cref{eq:chebyshev-lp}. Of all facets of $H_i$ tangential to this ball, choose the one with the largest associated dual variable in \cref{eq:chebyshev-lp} and take its normal vector. If already used, take the normal vector to the facet of $H_i$ with the largest volume.
\end{enumerate}

With all these methods, we also use a filtering procedure:
we discard vectors that have already been used.
For each method, it is clear how to generate another direction if the first was discarded: for example for the ``facets'' method we choose the normal of the facet with the second-largest volume if the first direction was discarded.
For the ``maximal-centre-then-facets'' methods we fall back on the ``facets'' method when all normals to facets tangantial to the Chebyshev ball have already been used.

In the filtering procedure, we employ an angle threshold $\theta$ such that a potential direction $d$ is discarded if it is within $\theta$ degrees of any previously used direction.
If the filter discards all possible directions, we reduce the angle threshold $\theta$ by 20\%.
This is repeated every time a method ``runs out of directions'', until $\theta$ falls below a pre-defined minimum angle $\theta_{min}$, at which point the whole algorithm is terminated.

Note that the approximation of $\Aeps$ can be parallelised effectively for any of the three direction generation methods, in the sense that multiple optimisations in different directions can be run in parallel.
For the latter two methods, this means that the convex hull $H_{i - P}$ must be used in the calculation of the $i$-th direction when there are $P$ parallel optimisations.
When $P$ is large, this means some of the generated direction could be slightly inferior (being generated with an older hull $H_{i - P}$).

\begin{figure}
	\centering
	\includegraphics{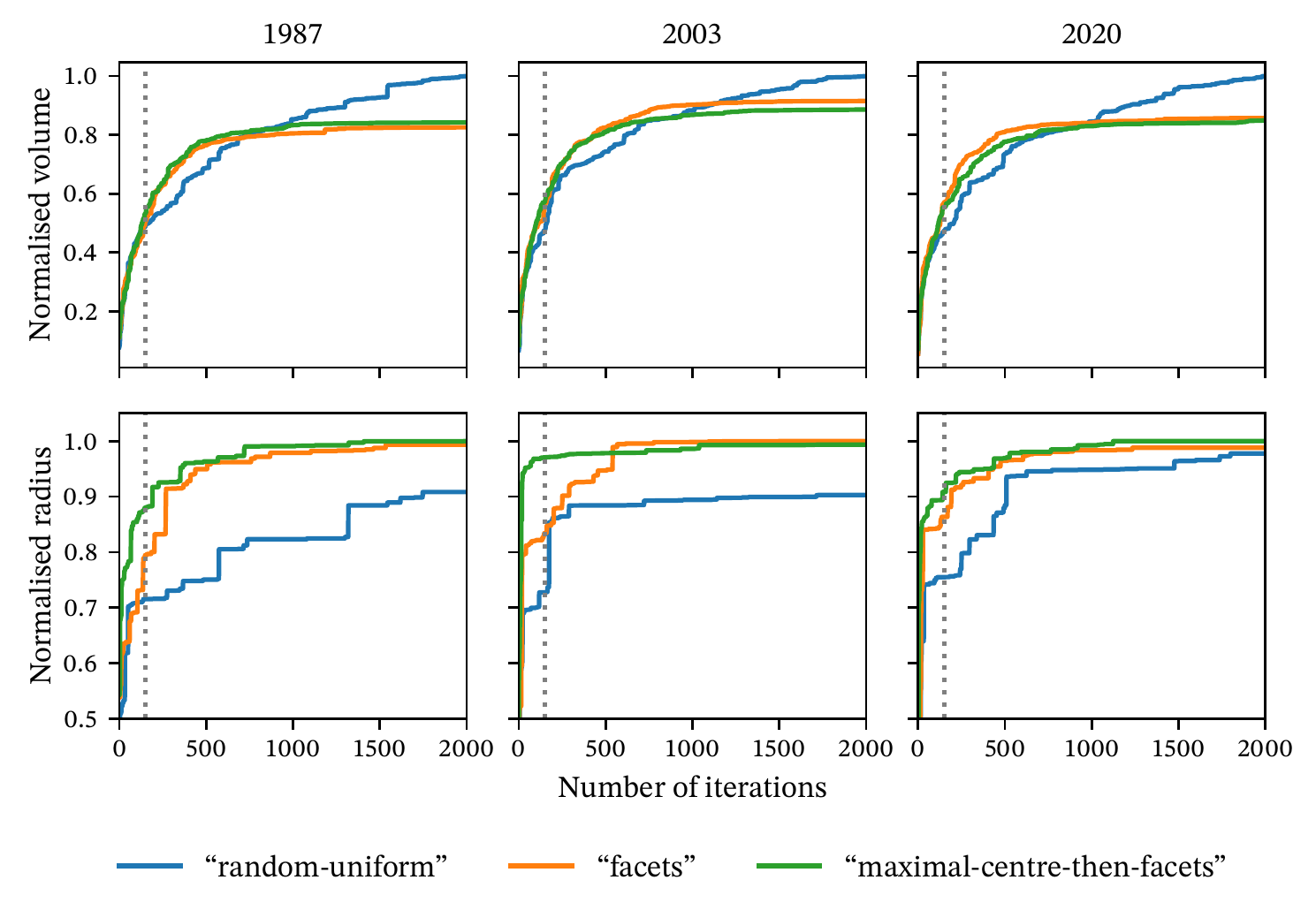}
	\caption{Performance of the different direction generation methods in terms of volume and Chebyshev radius convergence. We selected three years (among those, the years with the highest and lowest optimal costs) and approximated their near-optimal spaces with 2000 iterations. For each of the years the plotted volumes and radii have been normalised by the largest volume and radius obtained by any method for that year. The dotted lines mark 150 iterations.}
	\label{fig:direction-performance}
\end{figure}

We compare the performance of different methods in \cref{fig:direction-performance}.
The plots show that the three different direction generation methods have different characteristics, but also that their performance varies substantially between different spaces (different weather years).
Generally speaking, we see that the ``random-uniform'' method attains the largest volume in the long term, while the ``facets'' and ``maximal-centre-then-facets'' methods attain similar volumes and have a stronger performance initially.
In terms of radius, the ``random-uniform'' method performs worse, while the ``maximal-centre-then-facets' method converges the quickest initially.

We see that a large number of iterations is needed to converge in terms of volume, with the ``random-uniform'' attaining the highest volume, but still not converging after 2000 iterations.
Meanwhile, the other two methods based on facet normals make large strides initially, but display a false convergence below the actual volume after about 1000 iterations.
Convergence in terms of the radius is better (especially for the ``maximal-centre-then-facets'' method), but can be more erratic than the convergence of volume.
150 iterations with the ``maximal-centre-then-facest'' method were chosen as a compromise between accuracy and computational demand for this paper.

\begin{figure}
	\centering
	\includegraphics{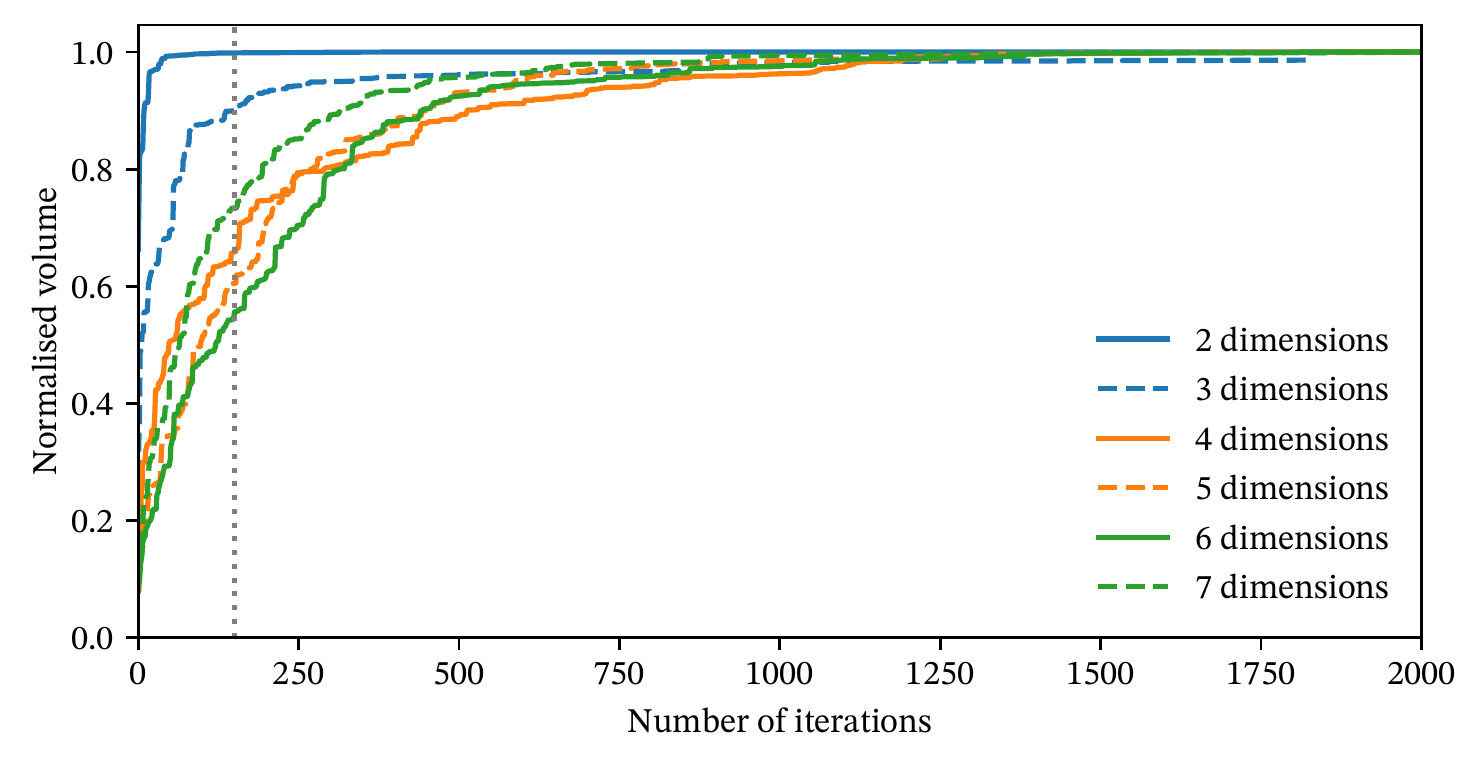}
	\caption{The convergence of volume in approximating the near-optimal space for the weather year 2020, reduced to different numbers of dimensions. The direction generation method is ``maximal-centre-then-facets''. The vertical dotted lines marks 150 iterations. The volume is normalised for each of the dimensions individually.}
	\label{fig:dim-iterations}
\end{figure}

In \cref{fig:dim-iterations} we compare convergence of volume between different numbers of dimensions $k$ of the spaces $\Aeps \subseteq \mathbb{R}^k$.
For this plot, we use the ``maximal-centre-then-facets'' direction generation method and approximated the near-optimal spaces for the weather year 2020, but use different dimension reduction maps $\sigma$.
The final approximated volumes are all normalised to 1.
We see that 3- and especially 2-dimensional spaces are quickly approximated, while the convergence is slower for higher dimensional spaces.
However, we do not find a significant difference in convergence between 4, 5, 6 or 7 dimensions.

The results in this section can be used to inform a termination criterion for \cref{alg:near-opt-approx}.
The simplest option is to terminate the algorithm after a fixed number of iterations.
Alternatively, the algorithm may be terminated after the volume or Chebyshev radius of $H_i$ have not changed more than $\delta$ percent between successive iterations for the last $N_{\text{conv}}$ iterations.
In this case, we advise that $N_{\text{conv}}$ be chosen as large as possible, since we can see from \cref{fig:direction-performance} that the convergence on volume and especially radius is often somewhat erratic.

\section{Data}
\label{sec:data-appendix}
\setcounter{figure}{0}
\renewcommand{\thefigure}{B.\arabic{figure}}
\subsection{Load data}
In this article we use load data based on two regressions that were trained on hourly country-level ENTSO-E data from 2010 to 2014 from \cite{entso-e-2022}\footnote{Due to inconsistencies for Swiss ENTSO-E data, we additionally used data from the Swiss transmission operator \cite{swissgrid-2022}.}.
The aim is to infer country-level synthetic load data for each weather year between 1980 and 2020, whose profiles relate to weather patterns but are otherwise directly comparable.
In other words, we disregard long-term changes in demand due to demographic and technological developments.

\begin{figure}
	\centering
	\includegraphics{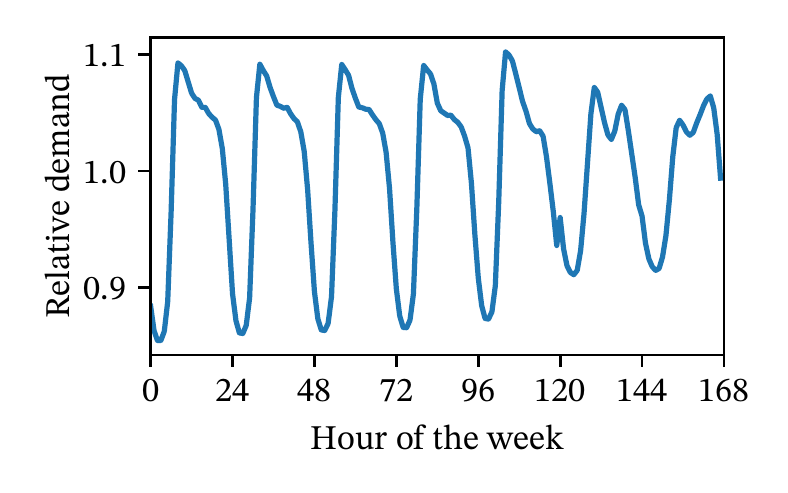}
	\caption{Weekly load profile for Norway, based on the regression described in \cref{eq:first-regression}.}
	\label{fig:weekly-profile}
\end{figure}

First we infer weekly load profiles (at an hourly resolution) for each country; for the purpose of this regression we treat holidays for each country as Sundays (using the Python package \texttt{python-holidays} \cite{montel-2022}).
Specifically, for each country $c$, we first divide the hourly demand values by the daily average value.
On these normalised values, we conduct a regression based on the following model formulation:
\begin{align}
	D_c^{\text{norm}}(t) = \alpha_c \cdot t + \alpha_{c,\: t \bmod 168},
	\label{eq:first-regression}
\end{align}
where $\alpha_c$ is an annual linear trend component, and the parameters $\alpha_{c,t \bmod 168}$ describe the weekly profile (see an example of this in \cref{fig:weekly-profile}).

\begin{figure}
	\centering
	\includegraphics{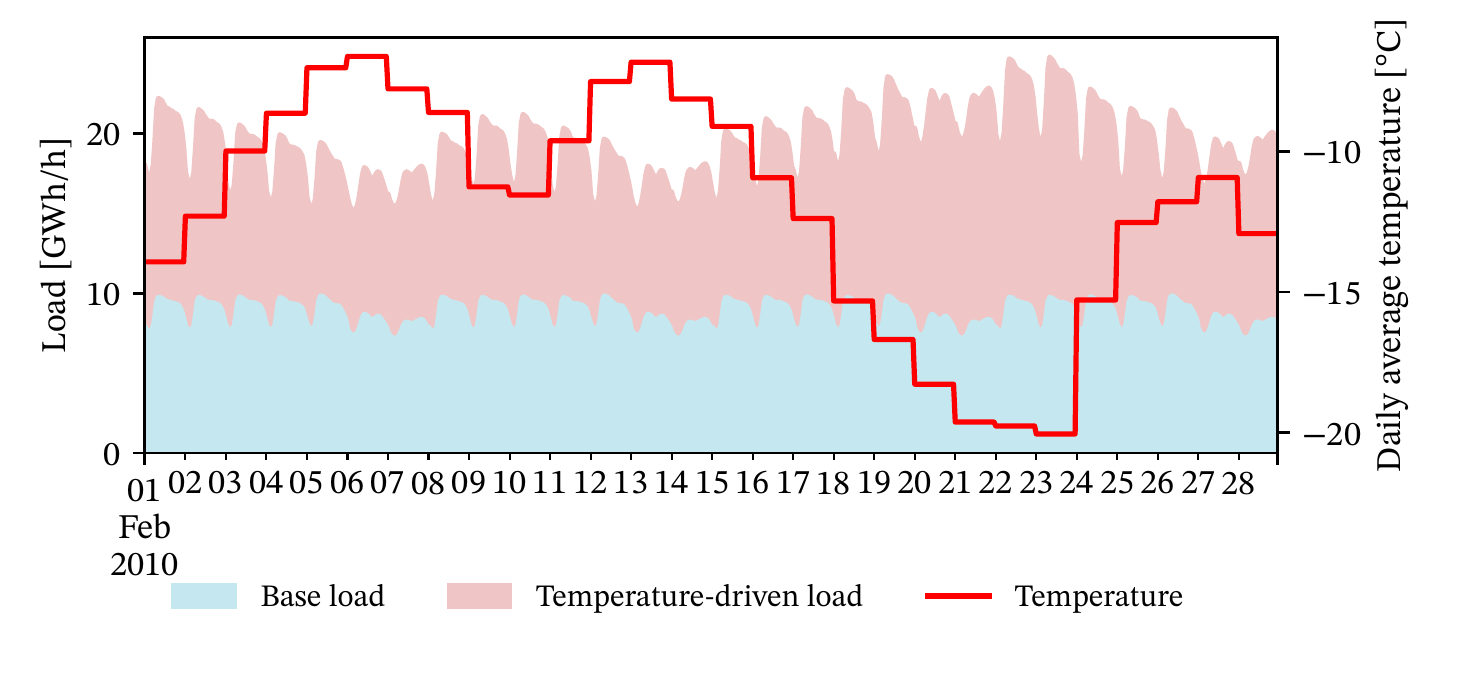}
	\caption{Load data (here for Norway in February 2010) split into base load and temperature-driven load. The underlying temperature data are daily, country-wide averages from ERA5 reanalysis. Base load and temperature-driven load are derived from the regression described in \cref{sec:data-appendix}.}
	\label{fig:load-data}
\end{figure}

Afterwards we use the concept of heating and cooling degree days (HDD and CDD resp.) as in \cite{benth-benth-ea-2008} to find temperature-independent daily demand and demand driven by heating or cooling demand.
For simplicity we only use one threshold for both HDD and CDD, which we set to be $15.5^\circ$C as \cite{spinoni-vogt-ea-2015} use for HDD:
\begin{align*}
	\text{CDD}_c(d) & = \max\{T_c(d) - 15.5, 0\}, \\
	\text{HDD}_c(d) & = \max\{15.5 - T_c(d), 0\},
\end{align*}
where $T_c(d)$ is daily average temperature in country $c$ during day $d$. Note that by using the daily averages we represent smoothing effects of thermal inertia on heating and cooling demand (compare \cref{fig:load-data}).
The country-wide temperatures are computed from ERA5 reanalysis data~\cite{era5-data} via the open-source tool Atlite \cite{hofmann-hampp-ea-2021}.
We now conduct a regression on the daily average load, $D_c(d)$, with dummy variables for each weekday (where national holidays are classified as Sundays), and exogenous variables given by heating degrees and cooling degrees:
\begin{align}
	D_c(d) = \beta_{c, \text{weekday}(d)} + \beta_c^{\text{cooling}} \cdot \text{CDD}_c(d) + \beta_c^{\text{heating}} \cdot \text{HDD}_c(d).
	\label{eq:second-regression}
\end{align}
For both regressions, we test the parameters for statistical significance; in some cases we thus set the trend parameter $\alpha_c$ and the cooling parameter $\beta_c^{\text{cooling}}$ to 0.

With these regressions and the temperatures for 1980--2020 from ERA5, we can compute the artificial load for each hour as follows:
\begin{equation} \label{eq:demand_art}
	\widetilde{D}_c(t) = \alpha_{c, t \bmod 168} \cdot \left[\alpha_c \cdot (t \bmod 8760) + \beta_{c, \text{weekday}(t)} + \beta^{\text{cooling}}_c \cdot \text{CDD}_c(t) + \beta^{\text{heating}}_c \cdot \text{HDD}_c(t)\right],
\end{equation}
where the index $c$ is over countries, $t$ is the time in hours,
$\alpha_{c,i}$ the regression parameter for the $i$-th hour of the week, $\alpha_c$ is the annual trend component, $\beta_{c,j}$ the regression parameter for weekday $j$, and $\beta^{\text{cooling}}_c$, $\beta^{\text{heating}}_c$ are the regression parameters for one degree of cooling/heating demand for country $c$.
We abuse notation slightly by writing $\text{CDD}_c(t)$ to mean $\text{CDD}_c(d)$ where $d$ is the day containing $t$ (and likewise for HDD).
For an example, see \cref{fig:load-data}.

We validate the regression on hourly ENTSO-E load data on the country level for 2015, and show it to be a good fit  --- more information about this can be found in the GitHub repository\footnote{\label{footnote:data-github}\url{https://github.com/aleks-g/multidecade-data/tree/v1.0}}.

In accordance with load projections for 2030 by the European commission\footnote{\url{https://ec.europa.eu/clima/document/download/ec1acac9-10fe-4eeb-915f-cad388990e0f_en}, Fig. 44 (accessed 23/06/2022)} we increased the demand in each country by 13\%.

\subsection{Hydropower data}
For hydropower data we follow the approach in PyPSA-Eur which uses ERA5 reanalysis data to generate inflow profiles (using Atlite) that are then scaled by historical country-level hydro generation data from the US Energy Information Administration (EIA) \cite{eia-hydro-data}.
We have extended the default dataset in PyPSA-Eur to cover the entire period of 1980 to 2020, and we have also normalised EIA's production data to EIA's capacity levels of 2020\textsuperscript{\ref{footnote:data-github}}:
\begin{align*}
	\text{gen}_c(y) = \text{nom\_gen}(y) \cdot \frac{\text{cap}_c(2020)}{\text{cap}_c(y)},
\end{align*}
where $\text{gen}_c(y)$ and $\text{nom\_gen}_c(y)$ are the normalised and reported hydropower generation for country $c$ in year $y$ respectively, and $\text{cap}_c(y)$ is the reported hydropower capacity for country $c$ in year $y$.

This allows a comparison throughout different years without any trends in infrastructure development.
To avoid anachronisms, we have distributed generation and capacities of former countries onto the current states (based on the first year of current borders, e.g. 1993 for Czechia and Slovakia, or the sum of West and East Germany)\textsuperscript{\ref{footnote:data-github}}.

The historical generation data ensure that the inflow profiles are scaled to reasonable values (see the general approach in \cite{schlachtberger-brown-ea-2017} and in particular Fig. 4 in \cite{liu-andresen-ea-2019}); we are interested in variability and not trends, therefore we want fixed capacities to have comparable weather years.

\begin{figure}
	\centering
	\includegraphics{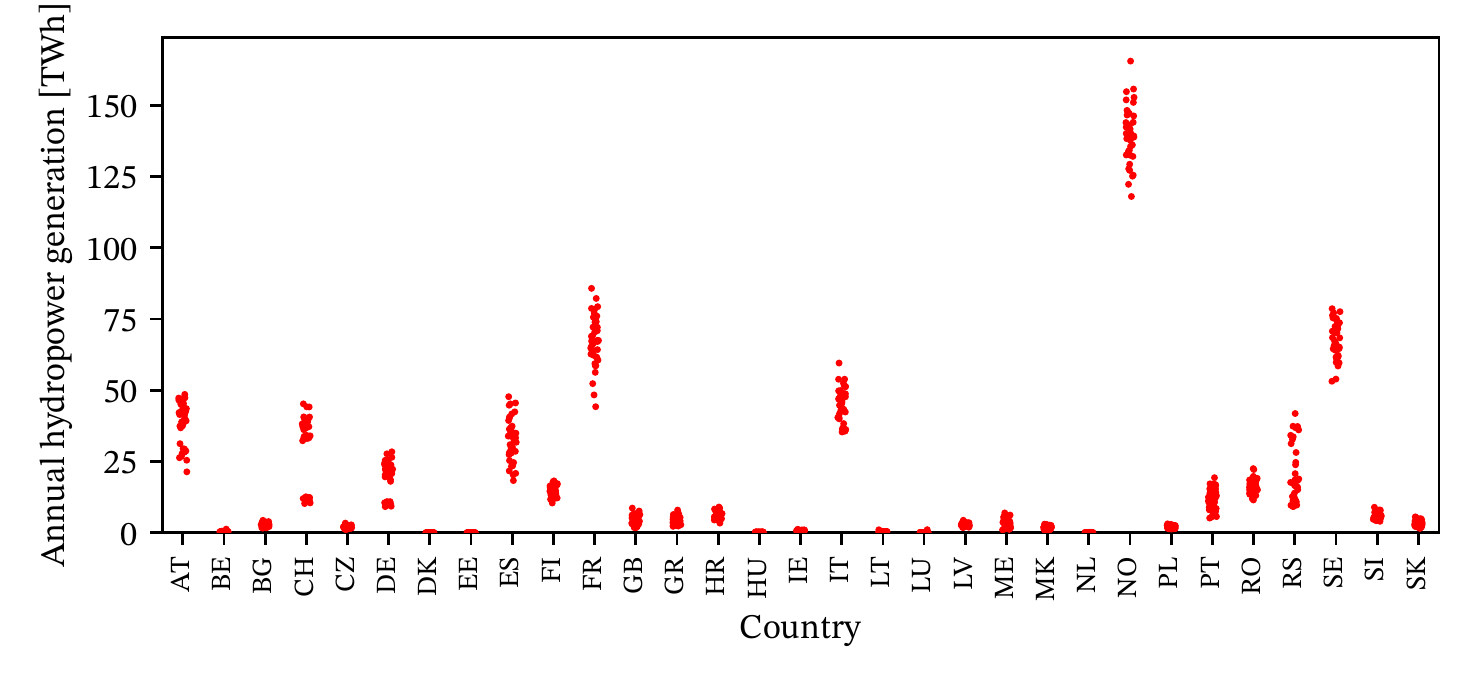}
	\caption{Variability of annual hydropower generation (1980--2020) based on EIA data \cite{eia-hydro-data}, normalised by the reported hydropower capacities from 2020.}
	\label{fig:hydro-data}
\end{figure}

\section{Additional results of a 100\% emission reduction}
\label{sec:emission-appendix}
\setcounter{figure}{0}
\renewcommand{\thefigure}{C.\arabic{figure}}

We showcase the impact of a 100\% emission reduction (as opposed to the 95\% reduction studied in the main text) on the results we present in this paper.
As the other assumptions remain unchanged and gas turbines are the only sources of \ch{CO2} emissions in our model, this corresponds to replacing gas as a generation technology with carbon-neutral alternatives.
The results sharpen the uncertainty that weather variability introduces to renewable power systems.
We investigate a different reduced near-optimal feasible space than before (replacing the gas investment dimension by two new dimensions, investment in battery and hydrogen storage).
This means that the Chebyshev centre under this new reduction is additionally robust to changes in investment in battery and hydrogen storage, which is not the case for the results in the main text.

\begin{figure}[tbh]
	\centering
	\includegraphics{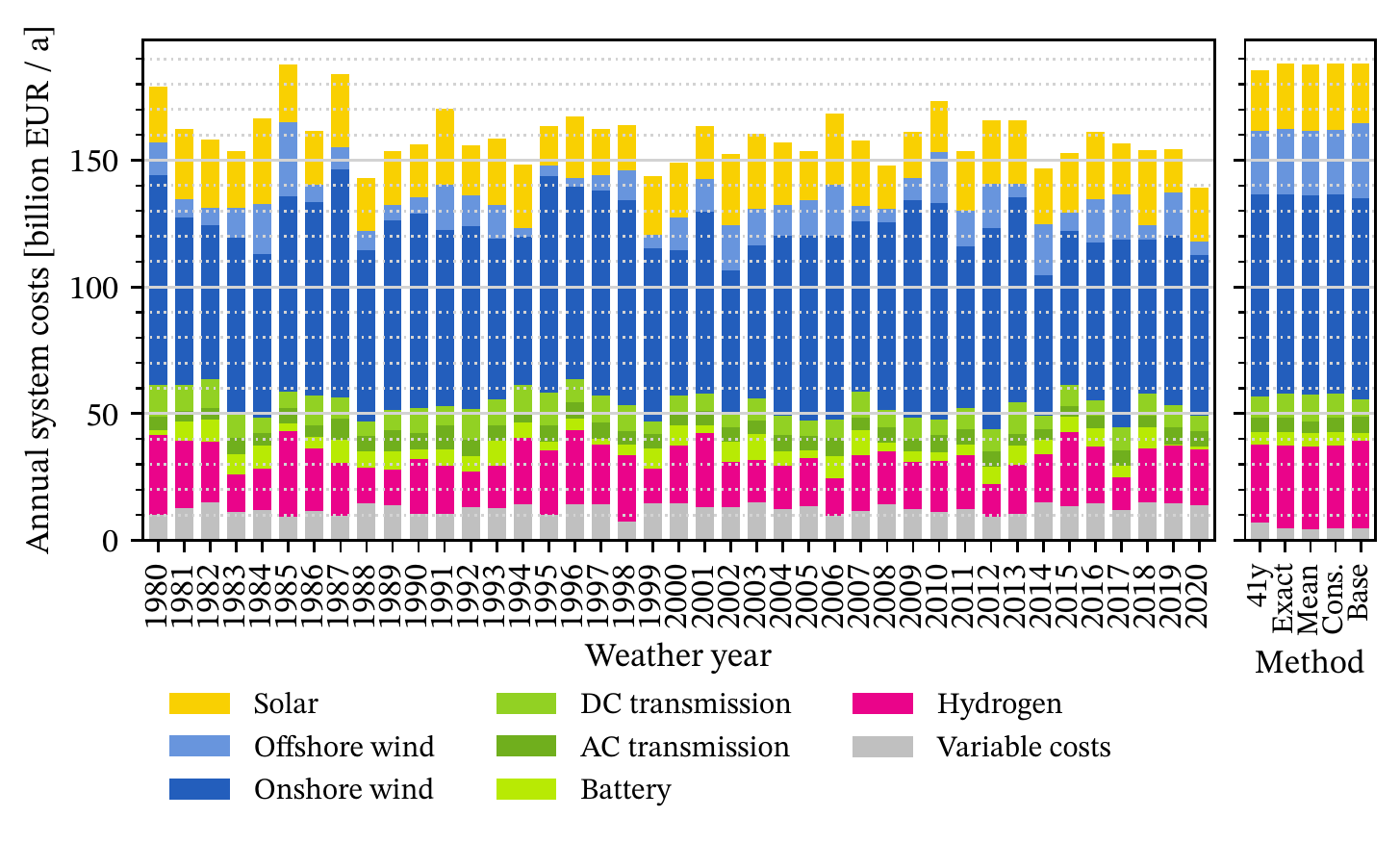}
	\caption{Comparison of cost-minimal designs under a 100\% emission reduction based on optimisations over individual years compared to the annualised costs of robust designs (exact, mean, conservative and baseline, respectively) and an optimisation with all 41 years (``41y'').}
	\label{fig:opt-costs-100_red}
\end{figure}

Without gas as a dispatchable generation technology, the costs of the now carbon-neutral power systems increase by 15 to 40 billion EUR/a, depending on the weather year.
The last percentage points of emission reductions are thus overproportionally costly (as also shown in \cite{neumann-brown-2021}), in particular in ``more difficult'' years.
Wind power in combination with storage technologies see an increase in investment (see \cref{fig:opt-costs-100_red}), most pronounced in the 1985, the year with the highest optimal costs.

\begin{figure}[tbh]
	\centering
	\includegraphics{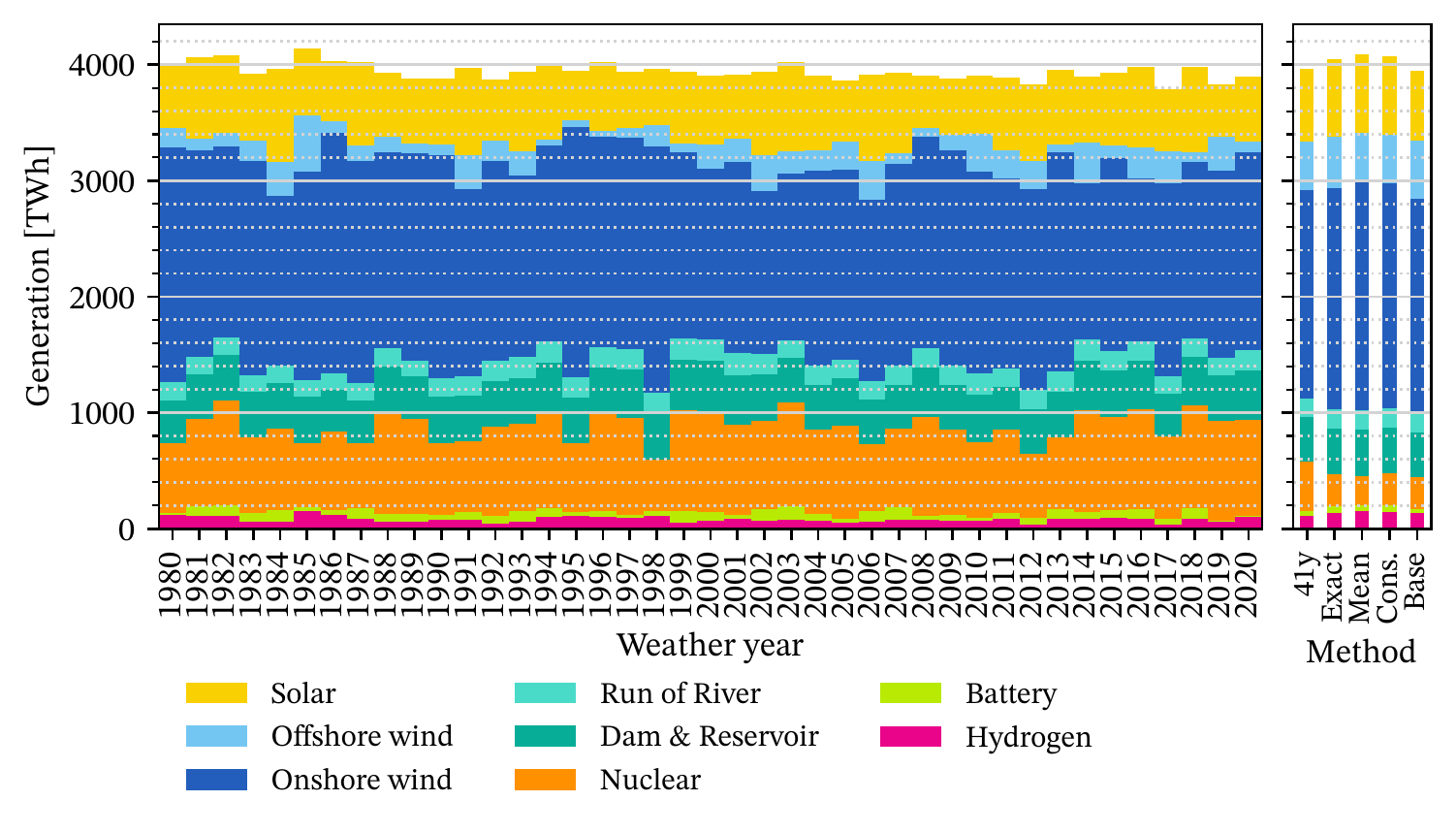}
	\caption{For 100\% emission reduction, annual average generation mixes for the optimal solutions for each single weather year, the optimal solution with all weather years (``41y'') and the robust allocations.}
	\label{fig:generation-100_red}
\end{figure}

\Cref{fig:generation-100_red} shows that the different optimal power systems are also mostly driven by (onshore) wind power, as in the 95\% reduction case.
Additionally, it depicts a strengthening of battery and hydrogen storage, which previously did not appear in the optimal solutions, although their shares of generation vary throughout the years.
Similarly, the optimal share of nuclear generation becomes more volatile and decreases not only in the robust allocations (which are characterised by higher renewable investment by design), but also in the 41-year optimisation.
As nuclear power has higher variable costs than renewable generators, the existing capacities are not fully used, indicating the competitiveness of newly installed renewable capacities.

\begin{figure}
	\centering
	\includegraphics{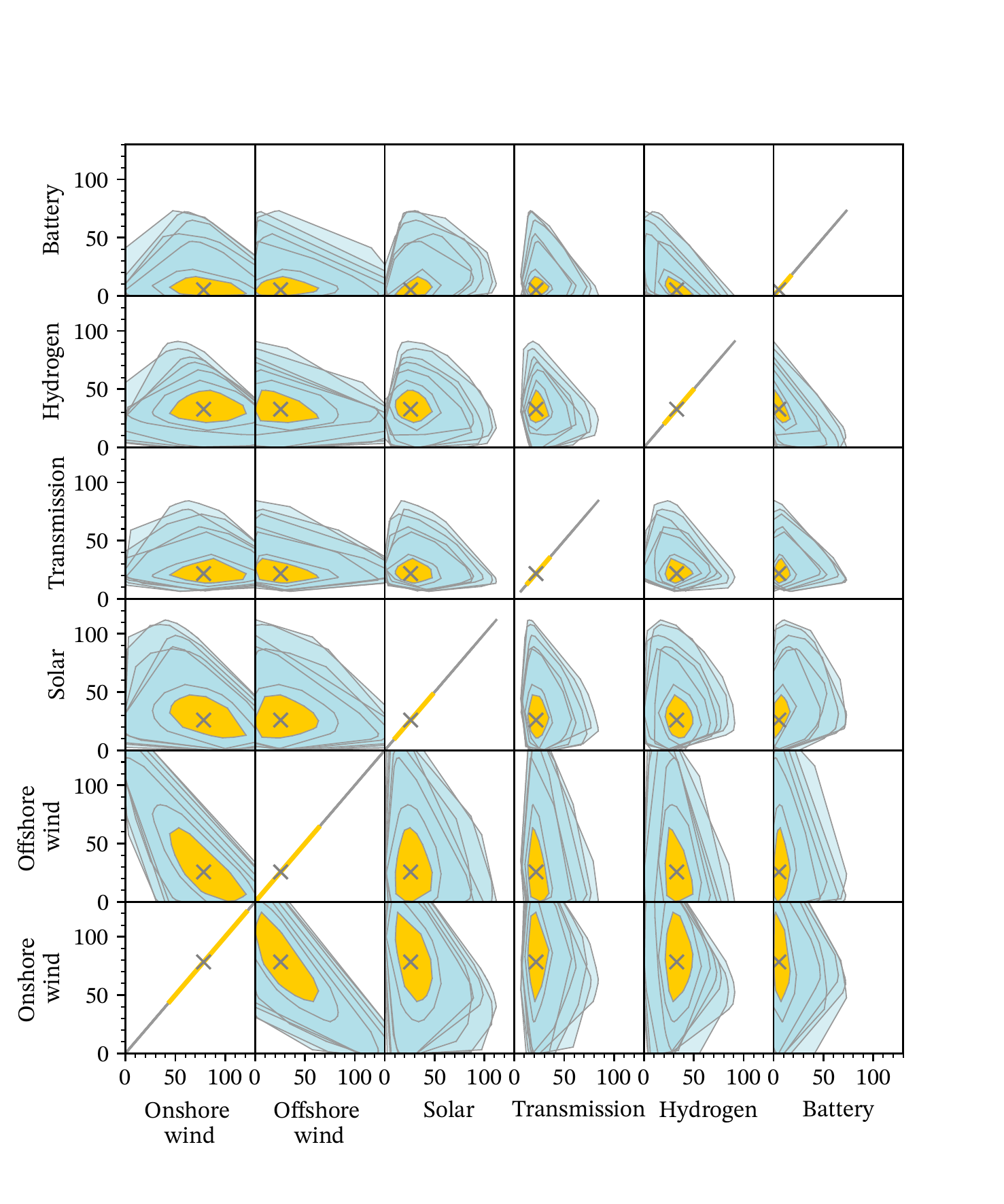}
	\caption{Projections of the near-optimal spaces for different weather years and their intersection under 100\% emission reduction. All values are annualised total investment costs per technology. For illustrative purposes, we only plot the near-optimal spaces for 6 out of 41 weather years (in different hues of blue). The intersection of all 41 near-optimal spaces is filled in yellow and the Chebyshev centre is marked with a cross.}
	\label{fig:near-opt-grid-100_red}
\end{figure}

\begin{figure}
	\centering
	\includegraphics{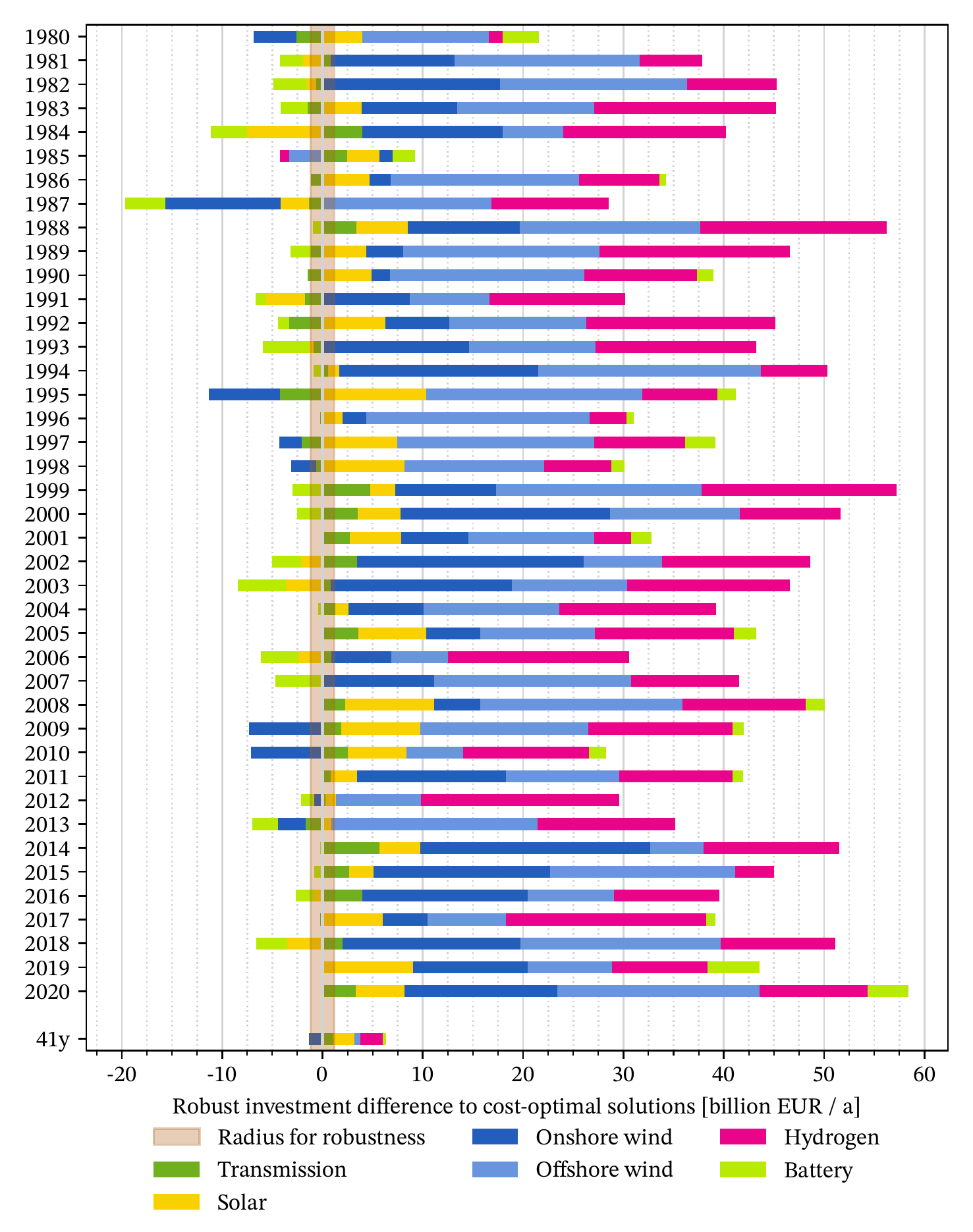}
	\caption{For 100\% emission reduction, a comparison of total investments for selected technologies in the optimal solutions for each weather year and the optimal solution with all weather years (``41y'') to the robust point $\ych$. Positive values mean greater investment in the given technology by the robust solution.}
	\label{fig:over-under-investment-100_red}
\end{figure}

As with the 95\% emission reduction, significant investment in both onshore wind power and solar power is necessary for a power system that can withstand weather variability.
\Cref{fig:opt-costs-100_red,fig:near-opt-grid-100_red} show additionally the need for investment in hydrogen storage and additional transmission capacities.
In comparison to optimal solutions for different weather years, additional investment in wind power and hydrogen strengthens robustness against weather variability (see \cref{fig:over-under-investment-100_red}).
The full decarbonisation increases the value of offshore wind power which did not feature prominently in the robust solution under 95\% emission reduction (see \cref{fig:over-under-investment}).
Finally, it should be noted that there is a significant amount of investment flexibility for onshore wind, offshore wind and solar among the robust solutions for all 41 weather years.

\end{document}